\documentclass[a4paper,12pt,reqno]{amsart}

\usepackage{amsmath}
\usepackage{amssymb}
\usepackage{amsfonts}
\usepackage{graphicx}
\usepackage{tabularx}
\usepackage{bm}
\usepackage[colorlinks]{hyperref}
\renewcommand\eqref[1]{(\ref{#1})} 
\newcommand{\dslash}{d\hspace{-0.4em}{ }^-\hspace{-0.2em}}

\newtheorem{theorem}{Theorem}[section]
\newtheorem{proposition}[theorem]{Proposition}
\newtheorem{lemma}[theorem]{Lemma}
\newtheorem{corollary}[theorem]{Corollary}
\theoremstyle{definition}
\newtheorem{definition}{Definition}[section]
\theoremstyle{remark}
\newtheorem*{remark}{Remark}
\graphicspath{ {images/} }

\newcommand{\eps}{\varepsilon}
\title[Hyperbolic systems with singular coefficients]{Hyperbolic systems with non-diagonalisable principal part and variable multiplicities, III: singular coefficients}
\author[Claudia Garetto]{Claudia Garetto}
\address{
	Claudia Garetto:
	\endgraf
	School of Mathematical Sciences
	\endgraf
	Queen Mary University of London
	\endgraf
	Mile End Road, London, E1 4NS
	\endgraf
	United Kingdom
	\endgraf
	{\it E-mail address} {\rm c.garetto@qmul.ac.uk}
}

\author[Bolys Sabitbek]{Bolys Sabitbek}
\address{
	Bolys Sabitbek:
	\endgraf
	School of Mathematical Sciences
	\endgraf
	Queen Mary University of London
	\endgraf
	Mile End Road, London, E1 4NS
	\endgraf
	United Kingdom
	\endgraf
	{\it E-mail address} {\rm b.sabitbek@qmul.ac.uk}
}

\thanks{The authors are supported by the
	EPSRC grant  EP/V005529/2}
\thanks{On behalf of all authors, the corresponding author states that there is no conflict of interest. This manuscript has no associated data.}

\date{}


\date{May 2022}

\begin{document}
	
	\maketitle
	
	\begin{abstract}
		In this paper we continue the analysis of non-diagonalisable hyperbolic systems initiated in \cite{GarJRuz, GarJRuz2}. Here we assume that the system has discontinuous coefficients or more in general distributional coefficients. Well-posedness is proven in the very weak sense for systems with singularities with respect to the space variable or the time variable. Consistency with the classical theory is proven in the case of smooth coefficients. 
	\end{abstract}
	
	\section{Introduction}
	This paper continues the analysis of hyperbolic systems with non-diagonalisable principal part and variable multiplicities started by Gramchev and Ruzhansky in \cite{GramRuz2013} for $2\times 2$ hyperbolic systems and continued by the first author, J\"ah and Ruzhansky in \cite{GarJRuz, GarJRuz2} for systems of arbitrary size. Here we assume that the system coefficients are singular, i.e., distributional in $x$ but still continuous with respect to $t$ or when smoothness is assumed in the space variable we allow singularities in time. Let us briefly summarise what is known about hyperbolic systems with non-diagonalisable principal part. Let  
	\begin{equation*}
	\begin{cases}
	D_t u=A(t,x,D_x)u+L(t,x,D_x)u+f(t,x),\qquad (t,x)\in[0,T]\times\mathbb{R}^n,\\
	u(0,x)=g(x),
	\end{cases}
	\end{equation*}
	where $A$ is an upper triangular $m\times m$ matrix with real eigenvalues and $L$ is the matrix of the lower order terms. Note that the eigenvalues might coincide and therefore we are dealing here in general with a hyperbolic system with multiplicities and non-diagonalisable principal part. It is also not restrictive to assume that the matrix $A$ is upper-triangular since conditions are given in \cite{GarJRuz} which allow the reduction into upper-triangular form of the system above. For a non exhausting overview on hyperbolic problems with multiplicities we refer the reader to \cite{B, ColKi:02, ColKi:02-2, CS, CDS, CJS2, CFDM:13, CFDM:13b, CFDM:15, dASpa:98, dAKi:05, Dencker, Garetto:15, GarRuz:1, GarRuz:2, GarRuz:3, GarRuz:4, GarRuz:6, GarRuz:7, GarJ, KY:06, KR2, OT:84, PP, Yagd}.
	
	Combining the results proven in \cite{GarJRuz, GarJRuz2,GramRuz2013} we have that the Cauchy problem above is well-posed in anisotropic Sobolev spaces of any order provided that the entries of $A=(a_{ij})$ are symbols in $C([0,T], S^1(\mathbb{R}^{2n}))$ and the entries of $L=(\ell_{ij})$ are symbols in $C([0,T], S^0(\mathbb{R}^{2n})$ fulfilling the following Levi type condition:
	\[
	\ell_{ij}\in C([0,T], S^{j-i}(\mathbb{R}^{2n})) 
	\]
	for $i>j$. In addition, when $A$ depends only on $x$ and assuming intersection of finite order at points of multiplicity one can obtain a representation formula for the solution which leads to results of propagation of singularities. See \cite{GarJRuz2} for more details. The results above heavily rely on the fact that the system coefficients are smooth in space and at least continuous on time. This is due to the fact that the system is solved making use of symbolic calculus and Fourier integral operators methods \cite{Hor:71}.
	
	When the system coefficients are less than continuous in time or less than smooth in space it might be difficult even to define a notion of solution. Indeed, as a consequence of Schwartz's impossibility result such system might fail to have a distributional solution \cite{HS:68, LO:91, Oberguggenberger:Bk-1992}. For this reason, we want to prove here well-posedness in a \emph{very weak sense}. Note that the desire to drop the regularity assumption of the system coefficients has physical motivation, since this kind of systems naturally arise in wave propagation into a multi-layered medium, like the Subsoil, or in conical refraction in crystals, see \cite{GarRuz:4} and references therein. 
	
	The notion of a \emph{very weak solution} has been introduced by the first author and Ruzhansky in \cite{GarRuz:4}. It is based on the idea of replacing the original problem with a regularised one, where the coefficients and initial data are regularised via convolution with suitable mollifiers. This generates a net of hyperbolic problems with smooth coefficients depending on the regularising parameter $\eps\in(0,1)$. The corresponding net of solutions $(u_\eps)_\eps$ is said to be a very weak solution if it satisfies some specific behaviour with respect to the parameter, that we will refer to as \emph{moderate}. Note that the notion of very weak solution is inspired by the theory of algebras of generalised functions \cite{GO:11a, GO:11b, GHO:06, GKOS:01, Oberguggenberger:Bk-1992} however, here we work with nets of functions rather than equivalence classes and we have full flexibility of regularising coefficients and initial data with different mollifiers and scales, basing our choice on the intrinsic properties of the problem we are studying. Moderateness will also be measured according to different norms or seminorms which are determined by the function space where we aim to prove well-posedness: $C^\infty$, Sobolev spaces, Gevrey classes, etc.
	
	In the sequel we list the main results obtained in this paper. 
	
	We begin by focusing on the Cauchy problem 
	\[
	D_t u=A(t,x,D_x)u+L(t,x,D_x)u+f(t,x),\qquad (t,x)\in[0,T]\times\mathbb{R}^n,
	\]
	with initial data $u(0,x)=g(x)$, where, 
	\begin{equation*}
	A(t,x,D_x) = \begin{bmatrix}
	\sum_{j=1}^n 	\lambda_{1j}(t,x) D_{x_j} & a_{12}(t,x,D_x) \\
	0 & \sum_{j=1}^n \lambda_{2j}(t,x)D_{x_j}
	\end{bmatrix},
	\end{equation*}
	and 
	\begin{equation*}
	L(t,x,D_x) = \begin{bmatrix}
	\ell_{11}(t,x) & \ell_{12}(t,x) \\
	\ell_{21}(t,x)\langle D_x\rangle^{-1} & \ell_{22}(t,x)
	\end{bmatrix}.
	\end{equation*}
	We assume that the matrices $A$ and $L$ are defined by pseudo-differential operators of order $1$ and $0$ respectively, with symbols which are not smooth in $x$. Namely, we assume that the system coefficients are continuous in $t$ but distributional with compact support in $x$. As a toy model one could think about a discontinuous function in $x$ or a Dirac distribution in $x$. After regularisation, i.e. convolution with a suitable mollifier
	\[
	\psi_\eps(x)=\omega(\eps)^{-n}\psi(\omega(\eps)^{-1}x),
	\]
	where $\omega(\eps)\to 0$ as $\eps\to 0$, we work on 
	\[
	\begin{cases}
	D_t u_{\varepsilon} = A_{\varepsilon} (t,x,D_x) u_{\varepsilon} + L_{\varepsilon} (t,x,D_x)u_{\eps} + f_{\varepsilon}(t,x), & \quad (t,x) \in [0,T]\times \mathbb{R}^n,\\
	u_{\varepsilon}(0,x) = g_{\varepsilon}(x), & \quad x \in \mathbb{R}^n,
	\end{cases}
	\]
	where $u_\eps(0,x) = [g_{1,\eps} (x), g_{2, \eps}(x)]^T$ and $f_\eps(t,x) = [f_{1,\eps}(t,x), f_{2, \eps}(t,x)]^T$, and the system matrices are
	\begin{equation*}
	A_{\varepsilon}(t,x,D_x) = \begin{bmatrix}
	\sum_{j=1}^n 	\lambda_{1j, \varepsilon}(t,x) D_{x_j} & a_{12, \varepsilon}(t,x,D_x) \\
	0 & \sum_{j=1}^n \lambda_{2j, \varepsilon}(t,x)D_{x_j}
	\end{bmatrix},
	\end{equation*}
	and 
	\begin{equation*}
	L_{\varepsilon}(t,x,D_x) = \begin{bmatrix}
	\ell_{11, \varepsilon}(t,x) & \ell_{12, \varepsilon}(t,x) \\
	\ell_{21, \varepsilon}(t,x)\langle D_x\rangle^{-1} & \ell_{22, \varepsilon}(t,x)
	\end{bmatrix}.
	\end{equation*}
	By the theory of Fourier integral operators \cite{GarJRuz, GarJRuz2} we know that the solution $u_\eps=[u_{1,\eps}, u_{2,\eps}]^T$ can be written explicitly via action of FIOs or integrated FIOs on the initial data, that we will denote here with $G^{0}_{i,\eps}$ and $G_{i,\eps}$, respectively, with $i=1,2$. The corresponding phase functions, $\phi_{i,\eps}$ can be also explicitly calculated by solving the Eikonal equation determined by the system eigenvalues. It turns out that to guarantee Sobolev mapping properties of the operators involved in our proof (see \cite{RuSu:06}) we need nets of symbols of order $0$ as well as some bound on the phase function, namely, 
	\begin{align*}
	|\partial_{\xi} \phi_{i, \eps}(t, x,\xi)-\partial_{\xi}\phi_{i, \eps}(t, y,\xi)| &\geq C(\eps)|x-y| \,\,\,\, \text{ for } \,\, x,y\in \mathbb{R}^n, \,\, \xi \in \mathbb{R}^n,  \\
	|\partial_y\phi_{i, \eps}(t, y,\xi)-\partial_y \phi_{i, \eps}(t, y,\eta)| &\geq C(\eps) |\xi - \eta| \,\,\,\, \text{ for } \,\,  y\in \mathbb{R}^n, \,\, \xi, \eta \in \mathbb{R}^n, \nonumber
	\end{align*}
	for $i=1,2$. In particular, it is crucial to study the $\eps$-behaviour of the composition operator 
	\[\mathcal{G}^0_{1, \eps} := G_{1, \eps} \circ (a_{12, \eps} + \ell_{12, \eps})\circ G_{2, \eps} \circ \ell_{21, \eps}\langle D_x\rangle^{-1}.
	\]
	Summarising, we identify the following set of hypotheses which guarantees the hyperbolicity of the system as well as its solvability in the very weak sense:
	\begin{itemize}
		\item[(H1)] the coefficients $(\lambda_{1j, \eps})_\eps$ and $(\lambda_{2j, \eps})_\eps$ are real valued and $C([0,T], C^\infty(\mathbb{R}^n))$-moderate for $j=1,\ldots,n$ with compact support in $x$;
		\item[(H2)]  $(\partial_k \lambda_{1j, \eps})_\eps$, $(\partial_k\lambda_{2j, \eps})_\eps$ and $(\ell_{ii, \eps})_{\eps}$ are logarithmic moderate for $k,j=1,\ldots,n$ and $i=1,2$;
		\item[(H3)] the nets of constants $(C(\eps))_\eps$ and $(C^{-1}(\eps))_\eps$ are moderate and \begin{equation*}
		|| \mathcal{G}^0_{1,\eps}  ||_{H^s \rightarrow H^s} = O(1).
		\end{equation*}
	\end{itemize}
	In detail, we prove that
	\begin{itemize}
		\item if (H1), (H2) and (H3) hold then the Cauchy problem above has a very weak solution (a solution $(u_\eps)_\eps$ which is $C([0,T], C^\infty(\mathbb{R}^n))$-moderate) and it is well-posed in the very weak sense (Theorem \ref{thm_general_statmenet}, Corollary \ref{cor_vws});
		\item if $A$ is diagonal and $L$ upper-triangular then (H1) and (H2) are sufficient to get very weak well-posedness (Corollary \ref{cor-1});
		\item if $A$ and $L$ are both upper-triangular then the Cauchy problem is well-posed in the very weak sense provided that (H1), (H2) hold and the nets of constants $(C(\eps))_\eps$ and $(C^{-1}(\eps))_\eps$ are moderate (Corollary \ref{cor-2});
		\item the last two results are easily extendable to matrices of size $m\times m$.
	\end{itemize}
	
	This paper is organised as follows: in Section 2 we collect some preliminaries on regularisation and very weak solutions for hyperbolic equations. The main results are proven in Section 3. In Section 4 we prove that very weak solutions recover the classical solutions in the limit as $\eps\to 0$ in the case of smooth coefficients. We also discuss some examples and we show how to adapt our methods to coefficients which are singular in time rather than in space. The paper ends with an Appendix devoted to the $L^2$-boundedness of Fourier integral operators \cite{RuSu:06}. In particular, we formulate precise bounds for FIOs depending on $\eps$ of the type encountered in the paper.
	
	\section{Preliminaries: regularisation and very weak solutions}
	In this section we recall the notion of very weak solution as introduced in \cite{GarRuz:4}. The main idea, which goes back to algebras of generalised functions, consists in regularising the equation coefficients via  a Friedrichs type mollifier $\psi$ ($\psi\in C^\infty_c(\mathbb{R}^n)$, $\psi\ge 0$ with $\int\psi=1)$. In other words, given a compactly supported distribution $u$ on $\mathbb{R}^n$, we generate a net of smooth functions $(u_\varepsilon)_\varepsilon$, where $\varepsilon\in (0,1]$, by writing 
	\[
	u_\varepsilon = u * \psi_{\omega(\varepsilon)},
	\]
	where 
	\[
	\psi_{\omega(\varepsilon)}(x) = \omega(\varepsilon)^{-n} \psi (x/\omega(\varepsilon)),
	\]
	and $\omega(\varepsilon)$ is a positive function converging to $0$ as $\varepsilon \rightarrow 0$ with $\omega(\varepsilon)\ge c\varepsilon^a$ for some $c,a>0$, uniformly with respect to $\varepsilon$. It is clear that $u_\varepsilon$ converges to $u$ as $\varepsilon\to 0$. In addition, the net $(u_\varepsilon)_\varepsilon$ fulfills some extra properties with respect to the parameter $\varepsilon$ that are illustrated in the proposition below. For a detailed proof of we refer the reader to \cite{GarRuz:4, GKOS:01, Oberguggenberger:Bk-1992}  and reference therein.
	\begin{proposition}
		\leavevmode
		\begin{itemize}
			\item[(i)] If $u  \in \mathcal{E'}(\mathbb{R}^n)$ then there exists $N \in \mathbb{N}_0$ and for all $\alpha \in \mathbb{N}_0^n$ there exists $c>0$ such that 
			\begin{equation*}
			|\partial^{\alpha} (u * \psi_{\omega(\varepsilon)}) (x)| \leq c \omega(\varepsilon)^{-N-|\alpha|},
			\end{equation*} 
			for all $x \in \mathbb{R}^n$ and $\varepsilon \in (0,1]$.
			\item[(ii)] If $f \in C_c^{\infty}(\mathbb{R}^n)$ then for all $\alpha \in \mathbb{N}_0^n$ there exists $c>0$ such that 
			\begin{equation*}
			|\partial^{\alpha}(f * \psi_{\omega(\varepsilon)})(x)|\leq c,
			\end{equation*}
			for all $x \in \mathbb{R}^n$ and $\varepsilon \in (0,1]$.
			\item[(iii)] If $f \in C_c^{\infty}(\mathbb{R}^n)$ and mollifier $\psi$ has all the moments vanishing, i.e., $\int_{\mathbb{R}^n} \psi(x)dx =1$ and $\int_{\mathbb{R}^n} x^{\alpha} \psi(x)dx =0$ for all multi-index $\alpha$ with $|\alpha|>1$, then for all $\alpha \in \mathbb{N}_0^{\alpha}$ and for all $q \in \mathbb{N}_0$ there exists $c>0$ such that 
			\begin{equation*}
			|\partial^{\alpha}(f*\psi_{\omega(\varepsilon)}(x)-f(x))|\leq c \omega(\varepsilon)^q, 
			\end{equation*} 
			for all $x \in \mathbb{R}^n$ and $\varepsilon \in (0,1]$.
		\end{itemize}
	\end{proposition}
	It follows that the nets obtained via convolution with a mollifier are \emph{moderate} or \emph{negligible} with respect to the seminorms of $C^\infty(\mathbb{R}^n)$ (see Definition \ref{def_mod_neg}). Analogously one can replace $C^\infty(\mathbb{R}^n)$ with $C^\infty([0,T]\times\mathbb{R}^n)$. In the sequel, $K \Subset \mathbb{R}^n$ stands for $K$ is a compact set in $\mathbb{R}^n$.  
	\begin{definition}
		\label{def_mod_neg}
		\leavevmode
		\begin{itemize}
			\item[(i)] A net $(v_{\varepsilon})_{\varepsilon}\in C^{\infty}(\mathbb{R}^n)^{(0,1]}$ is $C^{\infty}$-moderate if for all $K \Subset \mathbb{R}^n$ and for all $\alpha \in \mathbb{N}_0^n$ there exists $N \in \mathbb{N}_0$ and $c>0$ such that
			\begin{equation*}
			|\partial^{\alpha }v_{\varepsilon}(x)| \leq c \varepsilon^{-N},
			\end{equation*}
			uniformly in $x \in K$ and $\varepsilon \in (0,1]$.
			\item[(ii)] A net $(v_{\varepsilon})_{\varepsilon}\in C^{\infty}(\mathbb{R}^n)^{(0,1]}$ is $C^{\infty}$-negligible if for  all $K \Subset \mathbb{R}^n$, $\alpha \in \mathbb{N}_0^n$ and $q\in \mathbb{N}_0$ there exists $c>0$ such that
			\begin{equation*}
			|\partial^{\alpha} v_{\varepsilon}(x)| \leq c\varepsilon^q,
			\end{equation*}
			uniformly in $x \in K$ and $\varepsilon \in (0,1]$.
			\item[(iii)] A net $(v_{\varepsilon})_{\varepsilon} \in C^{\infty}([0,T] \times \mathbb{R}^n)^{(0,1]}$ is $C^{\infty}$-moderate if for all $K \Subset \mathbb{R}^n$ and for all $l \in \mathbb{N}_0$ and $\alpha \in \mathbb{N}_0^n$ there exists $N \in \mathbb{N}_0$ and $c>0$ such that
			\begin{equation*}
			|\partial^k_t\partial^{\alpha }v_{\varepsilon}(t,x)| \leq c \varepsilon^{-N},
			\end{equation*}
			uniformly in $t \in [0,T]$, $x \in K$ and $\varepsilon \in (0,1]$.
			\item[(iv)] A net $(v_{\varepsilon})_{\varepsilon} \in C^{\infty}([0,T] \times \mathbb{R}^n)^{(0,1]}$ is $C^{\infty}$-negligible if for  all $K \Subset \mathbb{R}^n$, $k \in \mathbb{N}_0$, $\alpha \in \mathbb{N}_0^n$ and $q\in \mathbb{N}_0$ there exists $c>0$ such that
			\begin{equation*}
			|\partial^k_t\partial^{\alpha} v_{\varepsilon}(t,x)| \leq c\varepsilon^q,
			\end{equation*}
			uniformly in $t \in [0,T]$, $x \in K$ and $\varepsilon \in (0,1]$.
			\item[(v)] A net $(v_{\varepsilon})_{\varepsilon} \in C([0,T], C^\infty(\mathbb{R}^n)^{(0,1]}$ is $C([0,T], C^\infty(\mathbb{R}^n))$-moderate if for all $K \Subset \mathbb{R}^n$ and for all $\alpha \in \mathbb{N}_0^n$ there exists $N \in \mathbb{N}_0$ and $c>0$ such that
			\begin{equation*}
			|\partial^{\alpha}_x v_{\varepsilon}(t,x)| \leq c \varepsilon^{-N},
			\end{equation*}
			uniformly in $t \in [0,T]$, $x \in K$ and $\varepsilon \in (0,1]$.
			\item[(vi)] A net $(v_{\varepsilon})_{\varepsilon} \in C^{\infty}([0,T] \times \mathbb{R}^n)^{(0,1]}$ is $C([0,T], C^\infty(\mathbb{R}^n))$-negligible if for  all $K \Subset \mathbb{R}^n$, $\alpha \in \mathbb{N}_0^n$ and $q\in \mathbb{N}_0$ there exists $c>0$ such that
			\begin{equation*}
			|\partial^{\alpha}_x v_{\varepsilon}(t,x)| \leq c\varepsilon^q,
			\end{equation*}
			uniformly in $t \in [0,T]$, $x \in K$ and $\varepsilon \in (0,1]$.
		\end{itemize}
	\end{definition} 
	\begin{remark}
		If we are dealing with nets of numbers then we will measure moderateness and negligibility in terms of the module in $\mathbb{C}$ or $\mathbb{R}$ and   we will simply talk of moderate and negligible nets. Analogously, we can state the definitions above for vectors by replacing $C^\infty(\mathbb{R}^n)$ with $C^\infty(\mathbb{R}^n)^m$. We will still talk of $C^\infty$-moderatenss and $C^\infty$-negligibility.
	\end{remark}
	\subsection{State of the art on hyperbolic equations with low regular coefficients}
	
	In this subsection we give a brief state of the art on hyperbolic equations with low regular coefficients, with particular focus on the very weak well-posedness of the corresponding Cauchy problem.
	
	We begin by recalling that the Cauchy problem for ordinary differential equations has been investigated already in \cite{GKOS:01} in the context of Colombeau algebras of generalised functions. The results obtained there (Theorem 1.5.2 \cite{GKOS:01}) can be rewritten at the net level as follows, for the Cauchy problem
	
	\begin{equation}
	\label{CPbook}
	\begin{split}
	\dot{x_\varepsilon}(t) &= F_{\varepsilon}(t,x_\varepsilon(t)),\\
	x_\varepsilon(t_0) &=x_{0,\varepsilon},
	\end{split}
	\end{equation}
	where $t\in\mathbb{R}$, $(F_\varepsilon)_\varepsilon$ is a net in $(C^\infty(\mathbb{R}^{1+n}))^n$ and $(x_{\varepsilon})_\varepsilon$ is a net of initial conditions. In the sequel we are making use of the standard norm in $\mathbb{R}^{1+n}$ that for the sake of simplicity we will continue to denote with $|\cdot|$.
	
	\begin{proposition}
		\label{theo_book}
		Let $(F_\varepsilon)_\varepsilon$ be a net in $(C^\infty(\mathbb{R}^{1+n}))^n$ such that 
		\begin{itemize}
			\item[(i)] for all $k\in\mathbb{N}_0$ and $\alpha\in\mathbb{N}_0^n$ there exists $N\in\mathbb{N}_0$ and $c>0$ such that 
			\[
			|\partial^k_t\partial^\alpha_x F_\varepsilon(t,x)|\le c\varepsilon^{-N}(1+|(t,x)|)^N,
			\]
			for all $(t,x)\in\mathbb{R}\times\mathbb{R}^n$ and $\varepsilon\in(0,1]$.
			\item[(ii)] there exists $c>0$ such that 
			\[
			\Vert \nabla_x F_\varepsilon(t,x)\Vert_{L^\infty(\mathbb{R}^{1+n})}\le c|\log \varepsilon|,
			\]
			for all $\varepsilon\in(0,1]$.
		\end{itemize}
		
		For all $\eps\in(0,1]$ the Cauchy problem \eqref{CPbook} has a unique solution $x_\eps\in C^\infty(\mathbb{R}^{1+n})^n$ such that if $(x_{0,\eps})_\eps$ is moderate then $(x_\eps)_\eps$ is $C^\infty(\mathbb{R})$-moderate as well. The solution $(x_\eps)_\eps$ is unique modulo, negligible perturbation, i.e., if we perturb the initial data with a negligible net (i.e., $x_{0,\eps}+n_\eps$, where $(n_\eps)_\eps$ is negligible) then the corresponding solution will differ from $(x_\eps)_\eps$ by a $C^\infty$-negligible net.

	\end{proposition}
	For the sake of simplicity we say that $(\nabla_x F_\eps)_\eps$ is logarithmic moderate. Note that one can reformulate the proposition above on a finite time interval $[0,T]$ and replace the polynomial bound in (i) with an assumption of compact support with respect to the $x$-variable.
	
	Proposition \ref{theo_book} can also be stated for nets that are only continuous with respect to $t$. This means to set $k=0$ and replace $C^\infty$ with $C([0,T], C^\infty(\mathbb{R}^n))$. Continuity with respect to $t\in[0,T]$ and smoothness with respect to $x$ will be the standard assumptions throughout this paper for our nets of functions.
	
	Let us now consider the scalar case, i.e., the Cauchy problem for first order hyperbolic equations. In detail, let 
	\begin{align}\label{eq-2.1}
	D_t w_\eps &= \sum_{j=1}^{n} a_{j, \varepsilon}(t,x) D_{j} w_{\varepsilon} + a_{0,\varepsilon}(t,x)w_{\varepsilon},\\
	w_{\varepsilon}(0,x) &= w_{0,\varepsilon}, \nonumber
	\end{align}
	where $D_j = D_{x_j}$, the coefficients $a_{j,\varepsilon}$ are real-valued for all $j=1,\ldots,n$. The equation coefficients are compactly supported distributions with respect to $x$ regularised by convolution with a mollifier $\psi$ that we assume non-negative and with integral $1$. In detail, 
	\begin{itemize}
		\item[-] $a_j \in C([0,T]\times \mathcal{E'}(\mathbb{R}^n))$ with $a_j \geq 0$ for all $j=1,\ldots, n$,
		\item[-] $a_0 \in C([0,T]\times \mathcal{E'}(\mathbb{R}^n))$
		\item[-] $w_0 \in \mathcal{E'}(\mathbb{R}^n)$.
	\end{itemize}
	and 
	\begin{itemize}
		\item[-] $a_{j, \varepsilon}(t,x) = (a_{j}(t, \cdot)*\psi_{\omega(\varepsilon)} )(x)$ for $j=1,\ldots,n$,
		\item[-] $a_{0, \varepsilon}(t,x) = (a_{0}(t, \cdot)*\psi_{\omega(\varepsilon)} )(x)$,
		\item[-] $w_0,\varepsilon(x) = w_0 * \psi_{\varepsilon}(x)$.
	\end{itemize}
	\begin{definition}
		\label{def_vws_1}
		We say that the Cauchy problem 
		\[
		\begin{split}
		D_t w &= \sum_{j=1}^{n} a_{j}(t,x) D_j w + a_{0}(t,x)w,\\
		w(0,x) &= w_{0}
		\end{split}
		\]
		admits a \emph{very weak solution} $(w_\eps)_\eps$ if the net $(w_\eps)_\eps$ solves the regularised problem \eqref{eq-2.1} and it is $C([0,T], C^\infty(\mathbb{R}^n))$-moderate.
	\end{definition}
	
	The following result, is the net-version of the well-posedness results obtained in \cite{GO:11a, GH:03, LO:91}. It shows that a transport equation with singular coefficients admits a very weak solution under suitable choice of the regularising scale. In the sequel we say that a net of function is logarithmic moderate if its $L^\infty([0,T]\times\mathbb{R}^n)$-norm is estimated by a logarithmic scale.
	\begin{proposition}
		\label{prop_GO}
		Let the coefficients $(a_{j,\varepsilon})_\eps$ be $C([0,T], C^\infty(\mathbb{R}^n))$-moderate for $j=1,\ldots,n$  with compact support in $x$.  Suppose that $(a_{j,\varepsilon})_\eps$ are real valued and $(\partial_k a_{j,\varepsilon})_\eps$ as well as $(a_{0,\varepsilon})_\eps$ are logarithmic moderate $(k,j=1,\ldots,n)$. Then,
		\begin{itemize}
			\item[(i)] the Cauchy problem above admits a very weak solution $(w_\eps)_\eps$ which is compactly supported with respect to $x$. 
			\item[(ii)] If $a_{0, \varepsilon} $ as well as $w_{0,\varepsilon}$ are real-valued then the solution $(w_\eps)_\eps$ is real valued as well.
			\item[(iii)] Negligible perturbations of the coefficients and initial data lead to a negligible perturbation in the very weak solution.
		\end{itemize}
	\end{proposition}
	In \cite{GO:11a} the authors investigated further the properties of the very weak solution $(w_\eps)_\eps$ and proved that it can be written as the action of a Fourier integral operator (FIO) on the initial data $(w_{0,\eps})_\eps$. We summarise here the main steps of their argument, mainly to solve the Eikonal equation to determine the phase function of the FIO and to solve a certain transport equation to determine its symbol. As expected, we will work with nets of FIO's, phase functions and symbols depending on the regularising parameter $\eps\in(0,1]$.
	
	\subsection{Nets of phase functions}
	We begin by constructing the phase function of the FIO making use of the characteristic curves associated to the principal part of the equation. This means to solve the linear Cauchy problem 
	\begin{equation}
	\begin{split}
	\label{Eikonal_CP}
	\partial_t \phi_\eps(t,x,\xi) &= \sum_{j=1}^n  a_{j,\varepsilon}(t,x)  \partial_j \phi_\eps(t,x,\xi)\\
	\phi_\eps(0,x,\xi) &= x \cdot \xi.
	\end{split}
	\end{equation}
	Under the assumptions on coefficients $(a_{j,\varepsilon})_\eps$ in Proposition \ref{prop_GO}, the Cauchy problem above is solved by a $C([0, T ], C^\infty(\mathbb{R}^n))$-moderate net $(\phi_\eps)_\eps$ which can be written as
	\begin{equation*}
	\phi_\eps(t,x,\xi) = \sum_{l=1}^n  \omega_{l,\eps}(t,x)\xi_l,
	\end{equation*}
	where $\omega_{l,\eps}$ with $l=1,\ldots,n$ are solutions of the Cauchy problem 
	\begin{align}\label{eq-3}
	\partial_t \omega_{l,\eps}(t,x) &= \sum_{j=1}^n a_{j,\varepsilon}(t,x) \partial_j\omega_{l,\eps}(t,x),\\
	\omega_{l,\eps}(0,x) &= x_l. \nonumber
	\end{align}
	
	We recall that the characteristic equations associated with the  Cauchy problem \eqref{eq-3} are given by
	\begin{align}
	\frac{d}{ds} \gamma_{l,\eps}(x,t;s) &= - a_{l,\varepsilon}(s,\gamma_{1,\eps}(x,t;s),\ldots,\gamma_{n,\eps}(x,t;s)),\\
	\gamma_{l,\eps}(x,t;t) &= x_l \quad \text{ for } \quad l=1,\ldots,n, \nonumber  
	\end{align}
	The solutions $\gamma_{1,\varepsilon},\ldots,\gamma_{n,\varepsilon}$ are the components of the characteristic curve 
	\[
	\gamma_\varepsilon =(\gamma_{1,\varepsilon},\ldots,\gamma_{n,\varepsilon})
	\]
	associated the differential operator $\sum_{j=1}^n a_{j,\varepsilon}(t,x) D_j$. The following proposition has been proven in \cite{GO:11a} in the context of generalised functions of Colombeau type. For the advantage of the reader we present here a detailed proof at the net level. 
	\begin{proposition}\label{prop-1}
		Let the coefficients $(a_{j,\varepsilon})_\eps$ be real valued $C([0, T ], C^\infty(\mathbb{R}^n))$-moderate, compactly supported in $x$ and let $(\partial_k a_{j,\varepsilon})_\eps$ be logarithmic moderate $(k,j=1,\ldots,n)$. Then, the solution $\omega_{l, \varepsilon}(t,x)$ of the Cauchy problem \eqref{eq-3} can be written as
		\begin{equation*}
		\omega_{l, \varepsilon}(x,t) = \gamma_{l, \varepsilon}(x,t; 0)= x_l + \int_{0}^t a_{l,\varepsilon}(\tau, \gamma_{1, \varepsilon}(x,t;\tau),\ldots,\gamma_{n, \varepsilon}(x,t;\tau))d\tau,
		\end{equation*}
		for $l=1,\ldots,n$.
	\end{proposition}
	
	\begin{proof}
		From Proposition \ref{prop_GO} it is easy to show that there exists a real valued solution $\omega_{l,\eps}$. It remains to prove that $\omega_{l,\eps}(t,x)= \gamma_{l,\eps}(x,t;0)$. This comes from the fact that $\omega_{l,\eps}$ is constant along the characteristic curves $\gamma_\eps(x,t;s)$ 
		\begin{align*}
		\frac{d}{ds} \omega_{l,\varepsilon}(s,\gamma_{1,\varepsilon}(x,t;s),\ldots,\gamma_{n,\varepsilon}(x,t;s))= 0.
		\end{align*}
		Hence, 
		\begin{align*}
		\omega_{l,\varepsilon}(t,\gamma_{1,\varepsilon}(x,t;t),\ldots,\gamma_{n,\varepsilon}(x,t;t)) = \omega_{l,\varepsilon}(0,\gamma_{1,\varepsilon}(x,t;0),\ldots,\gamma_{n,\varepsilon}(x,t;0))
		\end{align*}
		for each $t \in \mathbb{R}$. This implies $\omega_{l,\varepsilon}(t,x)=\gamma_{l,\varepsilon}(x,t;0)$. By solving directly the characteristic equations and making use of Proposition \ref{theo_book} we know that $\gamma_{l,\varepsilon}(x,t;0)$ is $C^{\infty}$-moderate and it has the following form 
		\begin{equation*}
		\gamma_{l, \varepsilon}(x,t;s) = x_l - \int_{t}^s a_{l,\varepsilon}(\tau, \gamma_{1, \varepsilon}(x,t;\tau),\ldots,\gamma_{n, \varepsilon}(x,t;\tau))d\tau,
		\end{equation*}
		then 
		\begin{equation*}
		\omega_{l,\varepsilon}(x,t) = \gamma_{l, \varepsilon}(x,t;0)= x_l + \int_{0}^t a_{l,\varepsilon}(\tau, \gamma_{1, \varepsilon}(x,t;\tau),\ldots,\gamma_{n, \varepsilon}(x,t;\tau))d\tau.
		\end{equation*}
		This completes the proof.
	\end{proof}
	As a straightforward consequence we obtain the following result for the Eikonal Cauchy problem \eqref{Eikonal_CP}

	\begin{corollary}\label{prop-2}
		If the coefficients $(a_{j,\varepsilon})_\eps$ are real valued $C([0, T ], C^\infty(\mathbb{R}^n))$-moderate with compact support in $x$ with $(\partial_k a_{j,\varepsilon})_\eps$ is logarithmic moderate $(k,j=1,\ldots,n)$ then the phase function 
		\begin{align*}
		\phi_\eps(t,x,\xi) &= \sum_{l=1}^n \omega_{l,\eps}(t,x)\xi_l = \sum_{l=1}^n \gamma_{l,\eps}(x,t;0)\xi_l, \\
		& =\sum_{l=1}^n x_l \xi_l+ \xi_l \int_{0}^t a_{l,\varepsilon}(\tau, \gamma_{1, \varepsilon}(x,t;\tau),\ldots,\gamma_{n, \varepsilon}(x,t;\tau))d\tau,
		\end{align*}
		solves the Eikonal Cauchy problem \eqref{Eikonal_CP}.
	\end{corollary}
	\subsection{The transport equation in this context}
	We now pass to solve the transport equation, in other words to solve the Cauchy problem \eqref{eq-2.1} with initial condition set to $1$. i.e.,  
	\begin{align}\label{eq-4}
	\begin{cases}
	D_t b_\eps(t,x) = \sum_{j=1}^n a_{j,\varepsilon}(t,x) D_j b_\eps(t,x) + a_{0, \varepsilon}(t,x)b_\eps(t,x), \\
	b_\eps(0,\cdot) = 1. 
	\end{cases}
	\end{align} 
	The characteristic equations of the Cauchy problem \eqref{eq-4} have the form 
	\begin{align*}
	\frac{d}{ds}\gamma_{l, \varepsilon}(x,t;s) &= - a_{l, \varepsilon}(s,\gamma_{1, \varepsilon}(x,t;s),\ldots, \gamma_{n, \varepsilon}(x,t;s)), \\
	\gamma_{l, \varepsilon}(x,t;t) &= 1.
	\end{align*} 
	and therefore
	\begin{equation*}
	\gamma_{l, \varepsilon}(x,t;s) = 1 - \int_{t}^{s} a_{l, \varepsilon} (\tau,\gamma_{1, \varepsilon}(x,t;\tau),\ldots, \gamma_{n, \varepsilon}(x,t;\tau) ) d\tau 
	\end{equation*}
	for $l = 1,\ldots, n$. By solving
	\begin{align*}
	\frac{d}{ds}z(s) &= a_{0, \varepsilon}(s, \gamma_{1, \varepsilon}(x,t;s),\ldots, \gamma_{n, \varepsilon}(x,t;s)) z(s), \\
	z(0) &= 1,
	\end{align*}
	where
	\begin{align*}
	z(s) = b_{\varepsilon}(s,\gamma_{1, \varepsilon}(x,t;s),\ldots, \gamma_{n, \varepsilon}(x,t;s)),
	\end{align*}
	we get
	\begin{equation}
	b_{\varepsilon}(t,x) = e^{i\int_{0}^{t}a_{0, \varepsilon}(\tau, \gamma_{1, \varepsilon}(x,t;\tau),\ldots, \gamma_{n, \varepsilon}(x,t; \tau))d\tau}.
	\end{equation}
	Using Proposition \ref{prop_GO} and integrating along the characteristics we have the following result.
	\begin{proposition}\label{prop-3}
		Let the coefficients $(a_{j,\varepsilon})_\eps$ be $C([0, T ], C^\infty(\mathbb{R}^n))$-moderate for $j=1,\ldots,n$  with compact support in $x$.  Suppose that $(a_{j,\varepsilon})_\eps$ are real valued and $(\partial_k a_{j,\varepsilon})_\eps$ as well as $(a_{0,\varepsilon})_\eps$ are logarithmic moderate $(k,j=1,\ldots,n)$, then the solution $(b_\eps)_\eps$ of the Cauchy problem \eqref{eq-4} is expressed by
		\begin{equation}
		b_\eps(t,x) = e^{i\int_{0}^{t}a_{0, \varepsilon}(s, \gamma_{1, \varepsilon}(t,x;s),\ldots, \gamma_{n, \varepsilon}(t,x;s))ds}
		\end{equation}
		and is $C([0, T ], C^\infty(\mathbb{R}^n))$-moderate.
	\end{proposition}
	We are now ready to write the solution of the Cauchy problem for first order hyperbolic equations in terms of a FIO formula. In other words we prove that the Cauchy problem in Definition \ref{def_vws_1} admits a very weak solution which is the action of a family of parameterised FIO's on the initial data. The results of the following subsection are extracted from \cite{GO:11a} where were originally stated in the context of generalised functions of Colombeau type.

	\subsection{Fourier Integral Operator formula}
	\begin{proposition}\label{prop-4}
		Let the coefficients $(a_{j,\varepsilon})_\eps$ be $C([0, T ], C^\infty(\mathbb{R}^n))$-moderate for $j=1,\ldots,n$  with compact support in $x$.  Suppose that $(a_{j,\varepsilon})_\eps$ are real valued and $(\partial_k a_{j,\varepsilon})_\eps$ as well as $(a_{0,\varepsilon})_\eps$ are logarithmic moderate $(k,j=1,\ldots,n)$. Then the solution $(w_{\eps})_\eps$ of the Cauchy problem \eqref{eq-2.1} is $C([0,T], C^\infty(\mathbb{R}^n))$-moderate and can be written as 
		\begin{equation}
		w_\eps(t,x) = F_{\phi_\eps}(b_\eps)(w_{0,\eps})(t,x):= \int_{\mathbb{R}^n} e^{i\phi_\eps(t, x, \xi)}b_\eps(t,x) \hat{w}_{0,\eps}(\xi) d\xi,
		\end{equation}
		where the phase function $\phi_\eps$ is defined by Corollary \ref{prop-2} and symbol $b_\eps$ by Proposition \ref{prop-3}.  
	\end{proposition}
	The solution of the inhomogeneous Cauchy problem 
	\begin{align}\label{eq-5}
	D_t w_\eps(t,x) &= \sum_{j=1}^n a_{j,\varepsilon}(t,x) D_j w_\eps(t,x) + a_{0,\varepsilon}(t,x)w_\eps(t,x) + f_{\varepsilon}(t,x),\\
	w_\eps(0,\cdot) &= w_{0,\eps}(x), \nonumber
	\end{align}
	where the net $(f_{\varepsilon})_\eps$ is $C([0, T ], C^\infty(\mathbb{R}^n))$-moderate with compact support in $x$, can also be expressed into FIO form. We begin by noting that Fourier integral operator $F_{\phi_\eps}(b_\eps) $ of Proposition \ref{prop-4} is given by   
	\begin{align*}
	F_{\phi_\eps}(b_\eps)(w_{0,\eps})(t,x) &= \int_{\mathbb{R}^n} e^{i\phi_\eps(t, x, \xi)}b_\eps(t,x) \hat{w}_{0,\eps}(\xi) d\xi \\
	& = b_\eps(t,x)\int_{\mathbb{R}^n} e^{i \gamma_{ \varepsilon}(x,t;0)\cdot \xi} \hat{w}_{0,\eps}(\xi) d\xi  \\
	&= b_\eps(t,x)w_{0,\eps}(\gamma_{\varepsilon}(x, t; 0)).
	\end{align*}
	It defines, for each $t \in \mathbb{R}$, a map 
	\begin{equation*}
	U_\eps(t) = F_{\phi_\eps}(b_\eps)(t) : w_{0,\eps} \rightarrow F_{\phi_\eps}(b_\eps)(w_{0,\eps})(t, \cdot), 
	\end{equation*}
	such that $U_\eps(0)=I$ and 
	\begin{equation*}
	U_\eps(t)^{-1}:v \rightarrow \frac{1}{b_\eps(t,\gamma_\eps(x,t;0))}v(\gamma_\eps(x,t;0)).
	\end{equation*}
	Hence, we obtain the following statement.
	
	\begin{theorem}
		\label{theo_FIO_inh}
		Let the coefficients $(a_{j,\varepsilon})_\eps$ be $C([0, T ], C^\infty(\mathbb{R}^n))$-moderate for $j=1,\ldots,n$  with compact support in $x$.  Suppose that $(a_{j,\varepsilon})_\eps$ are real valued and $(\partial_k a_{j,\varepsilon})_\eps$ as well as $(a_{0,\varepsilon})_\eps$ are logarithmic moderate $(k,j=1,\ldots,n)$. Then the solution $(w_{\eps})_\eps$ of the Cauchy problem \eqref{eq-5} is $C([0,T], C^\infty(\mathbb{R}^n))$-moderate and can be written as 
		\begin{equation}
		w_\eps(t,x) = F_{\phi}(b_\eps)(t) \left( w_{0, \eps} + i \int_{0}^{t} \frac{1}{b_\eps(\tau,\gamma_{\varepsilon}(\cdot,\tau;0))} f_{\varepsilon}(\tau, \gamma_{\varepsilon}(\cdot,\tau; 0))d\tau\right)(x),
		\end{equation}
		where the phase function $\phi_\eps$ is defined Corollary \ref{prop-2} and symbol $b_\eps$ by Proposition \ref{prop-3}. 
	\end{theorem}
	In the rest of the paper we will focus on hyperbolic systems with irregular coefficients. More precisely, we will formulate conditions on the principal part and the lower order terms matrix (Levi conditions) which will guarantee the existence of a very weak solution. Our analysis will also investigate different levels of singularity of the system coefficients and how they relate to the corresponding very weak solution.

	\section{Hyperbolic systems with irregular coefficients}
	We want to study the Cauchy problem for $2\times 2$ systems of the type
	\[
	D_t u=A(t,x,D_x)u+L(t,x,D_x)u+f(t,x),\qquad (t,x)\in[0,T]\times\mathbb{R}^n,
	\]
	with initial data $u(0,x)=g(x)$, where, 
	\begin{equation*}
	A(t,x,D_x) = \begin{bmatrix}
	\sum_{j=1}^n 	\lambda_{1j}(t,x) D_{x_j} & a_{12}(t,x,D_x) \\
	0 & \sum_{j=1}^n \lambda_{2j}(t,x)D_{x_j}
	\end{bmatrix},
	\end{equation*}
	and 
	\begin{equation*}
	L(t,x,D_x) = \begin{bmatrix}
	\ell_{11}(t,x) & \ell_{12}(t,x) \\
	\ell_{21}(t,x)\langle D_x\rangle^{-1} & \ell_{22}(t,x)
	\end{bmatrix},
	\end{equation*}
	are matrices of pseudo-differential operators of order 1 and 0, respectively, and $\lambda_{ij}\in\mathbb{R}$ for $i=1,2$ and $j=1,\ldots,n$. Differently from \cite{GarJRuz} we drop here the assumption of smoothness with respect to $x$, so the entries of the matrices above are assumed to be low regular in $x$, namely discontinuous or in general distributional. Passing to regularisation via convolution with a mollifier we are lead to study the regularised Cauchy problem
	\begin{equation}\label{Main-eq}
	\begin{cases}
	D_t u_{\varepsilon} = A_{\varepsilon} (t,x,D_x) u_{\varepsilon} + L_{\varepsilon} (t,x,D_x)u_{\eps} + f_{\varepsilon}(t,x), & \quad (t,x) \in [0,T]\times \mathbb{R}^n,\\
	u_{\varepsilon}(0,x) = g_{\varepsilon}(x), & \quad x \in \mathbb{R}^n,
	\end{cases}
	\end{equation}
	where $u_\eps(0,x) = [g_{1,\eps} (x), g_{2, \eps}(x)]^T$ and $f_\eps(t,x) = [f_{1,\eps}(t,x), f_{2, \eps}(t,x)]^T$, with  
	\begin{equation*}
	A_{\varepsilon}(t,x,D_x) = \begin{bmatrix}
	\sum_{j=1}^n 	\lambda_{1j, \varepsilon}(t,x) D_{x_j} & a_{12, \varepsilon}(t,x,D_x) \\
	0 & \sum_{j=1}^n \lambda_{2j, \varepsilon}(t,x)D_{x_j}
	\end{bmatrix},
	\end{equation*}
	and 
	\begin{equation*}
	L_{\varepsilon}(t,x,D_x) = \begin{bmatrix}
	\ell_{11, \varepsilon}(t,x) & \ell_{12, \varepsilon}(t,x) \\
	\ell_{21, \varepsilon}(t,x)\langle D_x\rangle^{-1} & \ell_{22, \varepsilon}(t,x)
	\end{bmatrix},
	\end{equation*}
	where $(\lambda_{1j, \eps})_\eps, (\lambda_{2j, \eps})_\eps \in  C([0,T], C^\infty(\mathbb{R}^n))^{(0,1]}$ for $j=1,\ldots, n$, are real-valued and $(a_{12, \eps})_\eps$ and $(\ell_{ik, \eps})_\eps$ are nets of symbols of order $1$ and $0$ respectively, i.e., $(a_{12, \eps})_\eps \in C([0,T], S^1_{0,1}(\mathbb{R}^n \times \mathbb{R}^n ))^{(0,1]}$ and $(\ell_{ik, \eps})_\eps \in C([0,T], S^0_{0,1}(\mathbb{R}^n \times \mathbb{R}^n ))^{(0,1]}$ for $i,k=1,2$.
	
	We will prove that the original Cauchy problem admits a very weak solution by carefully analysing the Sobolev well-posedness of the regularised Cauchy problem \eqref{Main-eq}. This will require some preliminary work inspired by \cite{GarJRuz} and the notions of Sobolev moderateness and negligibility.

	\subsection{Preliminaries: Fourier integral operators and Sobolev norms}\label{Aux_rem}
	For each $\lambda_{ij, \eps}$, $i=1,2$, $j=1,\ldots,n$ we will be denoting by $G_{i, \eps}^0 w_{0, \varepsilon}$ the solution to  
	\begin{align}\label{eq2.1}
	\begin{cases}
	D_t w_{\varepsilon} =  \sum_{j=1}^n \lambda_{ij, \varepsilon}(t,x) D_{x_j} w_{\varepsilon} + \ell_{ii, \varepsilon}(t,x, D_x)w_{\varepsilon},\\
	w_{\varepsilon}(0,x) = w_{0, \varepsilon}(x),
	\end{cases}
	\end{align}
	and by $G_{i, \eps} f_{\varepsilon}$ the solution to 
	\begin{align*}
	\begin{cases}
	D_t w_{\varepsilon} = \sum_{j=1}^n \lambda_{ij, \varepsilon}(t,x) D_{x_j} w_{\varepsilon} + \ell_{ii, \varepsilon}(t,x, D_x)w_{\varepsilon} + f_{\varepsilon}(t,x),\\
	w_{\varepsilon}(0,x) = 0. \nonumber
	\end{cases}
	\end{align*}
	These are Fourier integral operators where the phase function $\phi_{i, \eps}$ is defined by Corollary \ref{prop-2} and symbol
	\[
	b_{i, \eps}(t,x)=e^{i\int_{0}^{t}\ell_{ii,\eps}(s, \gamma_{1, \varepsilon}(t,x;s),\ldots, \gamma_{n, \varepsilon}(t,x;s))ds},
	\]
	by Proposition \ref{prop-3} for $i=1,2$. Namely, for
	\[
	\phi_{i,\eps}(t,x,\xi) =  x\xi+ \sum_{j=1}^n \xi_j \int_{0}^t \lambda_{ij,\varepsilon}(\tau, \gamma_{\eps}(x,t;\tau))d\tau=\gamma_{i,\eps}(x,t;0)\xi,
	\]
	the operators $G^0_{i, \eps}$ and $G_{i, \eps}$, $i=1,2$, are given by
	\begin{equation*}
	\begin{split}
	G^0_{i, \eps} w_{0, \varepsilon}(t,x) = F_{\phi_{i, \varepsilon}}(b_{i, \varepsilon})(w_{0, \varepsilon})(t,x)&:= \int_{\mathbb{R}^n} e^{{\rm i}\phi_{i, \varepsilon}(t, x, \xi)}b_{i, \varepsilon}(t,x) \hat{w}_{0, \varepsilon}(\xi) d\xi,\\
	&=b_{i,\eps}(t,x)w_{0,\eps}(\gamma_{i,\eps}(x,t;0))
	\end{split}
	\end{equation*}
	and
	\begin{align*}
	&G_{i, \eps} f_{\varepsilon}(t,x) =  F_{\phi_{i, \varepsilon}}(b_{i, \varepsilon})(t) \left( i \int_{0}^{t} \frac{1}{b_{i, \varepsilon} (s,\gamma_{i, \varepsilon}(\cdot, s;0))} f_{\varepsilon}(s, \gamma_{i, \varepsilon}(\cdot,s; 0))ds\right)(x)\\
	&=b_{i,\eps}(t,x) \left( i \int_{0}^{t} \frac{1}{b_{i, \varepsilon} (s,\gamma_{i, \varepsilon}(\gamma_{i, \varepsilon}(x,t; 0),s;0)} f_{\varepsilon}(s, \gamma_{i, \varepsilon}(\gamma_{i, \varepsilon}(x,t; 0),s; 0)))ds\right)\\
	&= i \int_{0}^{t}\int_{\mathbb{R}^{2n}}{\rm e}^{i\gamma_{i, \varepsilon}(\gamma_{i, \varepsilon}(x,t; 0),s; 0)\xi-iy\xi} \frac{1}{b_{i, \varepsilon} (s,\gamma_{i, \varepsilon}(\gamma_{i, \varepsilon}(x,t; 0),s;0))}b_{i,\eps}(t,x)f_{\varepsilon}(s,y)dy \dslash \xi ds\\
	&=\int_{0}^{t}\int_{\mathbb{R}^{2n}}{\rm e}^{i\gamma_{i, \varepsilon}(\gamma_{i, \varepsilon}(x,t; 0),s; 0)\xi-iy\xi} A_{i,\varepsilon}(s,t,x) f_{\varepsilon}(s,y)dy \dslash \xi ds.
	\end{align*}
	
	Combining Proposition \ref{prop-4} with Theorem \ref{theo_FIO_inh} we easily obtain the following proposition which describes the mapping properties of the operators $G^0_{i,\eps}$ and $G_{i,\eps}$.
	\begin{proposition}
		\label{prop_mapping_properties}
		
		Assume that
		\begin{itemize}
			\item[(H1)] The coefficients $(\lambda_{ij, \eps})_\eps$ are $C([0,T], C^{\infty}(\mathbb{R}^n))$-moderate for $i=1,2$ and $j=1,\ldots, n$ with compact support in $x$.
			\item[(H2)] The coefficients $(\lambda_{ij, \eps})_\eps$ are real valued and $(\partial_k\lambda_{ij, \eps})_\eps $ as well as $(\ell_{ii,\eps})_\eps$ are logarithmic moderate for $k,j=1,\ldots,n$ and $i=1,2$. 
		\end{itemize}
		The operators $G^0_{i,\eps}$ and $G_{i,\eps}$
		map $C([0,T], C^{\infty}(\mathbb{R}^n))$-moderate nets with compact support in $x$ into nets of the same type. The same holds with moderate replaced by negligible. 
	\end{proposition}
	We recall that pseudo-differential operators are a special kind of Fourier integral operators so in the sequel we will also deal with nets of pseudo-differential operators which will therefore have the same mapping properties of $G^0_{i,\eps}$. More precisely, we will work with the following nets of symbols.
	\begin{definition}
		\leavevmode
		\begin{itemize}
			\item[(i)]      A net $(a_\eps)_\eps \in C^{\infty}(\mathbb{R}^n \times \mathbb{R}^n )^{(0,1]}$ is $S^m_{0,1}(\mathbb{R}^n \times \mathbb{R}^n )$-moderate if for all $K \Subset \mathbb{R}^n$ and for all $\alpha, \beta \in \mathbb{N}_0^n$ and there exists $N \in \mathbb{N}_0$ and $c>0$ such that
			\begin{equation*}
			\langle \xi \rangle^{|\beta| - m}| \partial^{\alpha}_x \partial^{\beta}_{\xi} a_{\eps}(x,\xi)| \leq  c \eps^{-N},
			\end{equation*}
			uniformly in $\xi \in \mathbb{R}^n$, $x \in K$ and $\eps \in (0,1]$.
			\item[(ii)]     A net $(a_\eps)_\eps \in C^{\infty}(\mathbb{R}^n \times \mathbb{R}^n )^{(0,1]}$ is $S^m_{0,1}(\mathbb{R}^n \times \mathbb{R}^n )$-negligible if for all $K \Subset \mathbb{R}^n$ and for all $\alpha, \beta \in \mathbb{N}_0^n$ and there exists $q \in \mathbb{N}_0$ and $c>0$ such that
			\begin{equation*}
			\langle \xi \rangle^{|\beta| - m}| \partial^{\alpha}_x \partial^{\beta}_{\xi} a_{\eps}(x,\xi)| \leq  c \eps^{q},
			\end{equation*}
			uniformly in $\xi \in \mathbb{R}^n$, $x \in K$ and $\eps \in (0,1]$.
		\end{itemize}
	\end{definition}
	
	
	Finally, we denote by $C([0,T], S^m_{0,1}(\mathbb{R}^n \times \mathbb{R}^n ))^{(0,1]}$ the space of net of all symbols $(a_\eps(t,x,\xi))_\eps \in S^m_{0,1}(\mathbb{R}^n \times \mathbb{R}^n )^{(0,1]}$ which are continuous with respect to $t$ and in analogy with the definition above we can define moderate and negligible nets in $C([0,T], S^m_{0,1}(\mathbb{R}^n \times \mathbb{R}^n ))^{(0,1]}$ by employing uniform estimates with respect to $t\in[0,T]$.
	
	We can now go back to our regularised Cauchy problem. As explained in detail in \cite{GarJRuz} the solution of the Cauchy problem \eqref{Main-eq} can be written as
	\begin{eqnarray}
	\label{eq:u21} u_{1,\varepsilon} &=& U^0_{1,\varepsilon} + G_{1,\eps}((a_{12,\eps}+\ell_{12,\eps})u_{2,\eps}), \\
	\label{eq:u22} u_{2,\eps} &=& U^0_{2,\eps} + G_{2,\eps}(\ell_{21,\eps}\langle D_x \rangle^{-1} u_{1,\eps}), 
	\end{eqnarray} where \begin{equation} \label{eq:aux52}
	U_{i,\eps}^0 = G_{i,\eps}^0 g_{i,\eps} + G_{i,\eps}(f_{i,\eps}), \quad i=1,2.
	\end{equation} 
	We want to prove that a smooth solution exists by proving that it belongs to every Sobolev space $H^s$. This means to make use of the Sobolev mapping properties of all the operators involved above.
	
	Plugging \eqref{eq:u22} in \eqref{eq:u21}, we obtain 
	\begin{equation} \label{eq:u1Fin2x2}
	u_{1,\eps} = \tilde{U}^0_{1,\eps} + G_{1,\eps}(a_{12,\eps}G_{2,\eps}(\ell_{21,\eps}\langle D_x \rangle^{-1} u_{1,\eps})) + G_{1,\eps}(\ell_{12,\eps}G_{2,\eps}(\ell_{21,\eps}\langle D_x \rangle^{-1}u_{1,\eps})),
	\end{equation} 
	where 
	\begin{equation} \label{eq:aux51}
	\tilde{U}^0_{1,\eps} = G_{1,\eps}^0 g_{1,\eps} + G_{1,\eps}(f_{1,\eps}) + G_{1,\eps}((a_{12,\eps}+\ell_{12,\eps})U^0_{2,\eps}).
	\end{equation}
	
	By employing Proposition \ref{Prop-L2} in the Appendix we have that 
	\begin{equation}
	\label{est_L_2}
	\begin{split}
	\Vert G_{i,\eps}v\Vert_{H^s}&\le TC'(\eps)\sup_{t,s\in[0,T]}\sup_{|\beta|\leq 2n+1} || \partial_x^{\beta} A_{i, \varepsilon}(s,t,x)||_{L^{\infty}(\mathbb{R}^n_x)}\Vert v\Vert_{H^s},\\
	\Vert G^0_{i,\eps}v\Vert_{H^s} &\le C'(\eps)\sup_{t\in[0,T]}\sup_{|\beta|\leq 2n+1} || \partial_x^{\beta} b_{i, \varepsilon}(t,x)||_{L^{\infty}(\mathbb{R}^n_x)}\Vert v\Vert_{H^s},
	\end{split}
	\end{equation}
	where $C'(\eps)$ depends on the positive nets $C(\eps)$, $C_{\alpha,i}(\eps)$ and $C_{\beta,i}(\eps)$ appearing 
	below:
	\begin{align*}
	|\partial_{\xi} \phi_{i, \eps}(t, x,\xi)-\partial_{\xi}\phi_{i, \eps}(t, y,\xi)| &\geq C(\eps)|x-y| \,\,\,\, \text{ for } \,\, x,y\in \mathbb{R}^n, \,\, \xi \in \mathbb{R}^n,  \\
	|\partial_y\phi_{i, \eps}(t, y,\xi)-\partial_y \phi_{i, \eps}(t, y,\eta)| &\geq C(\eps) |\xi - \eta| \,\,\,\, \text{ for } \,\,  y\in \mathbb{R}^n, \,\, \xi, \eta \in \mathbb{R}^n, \nonumber
	\end{align*}
	and	
	\begin{equation*}
	\begin{split}
	|\partial_y^{\alpha}\partial_{\xi}\phi_{i, \eps}(t, y,\xi)| &\leq C_{\alpha, i}(\eps),\\	|\partial_y\partial_{\xi}^{\beta}\phi_{i, \eps}(t, y,\xi)| &\leq C_{\beta, i}(\eps),
	\end{split}
	\end{equation*}
	for all $t\in[0,T]$, for $1\leq |\alpha|, |\beta|\leq 2n+2$ and $i=1,2$.

	In the sequel, we introduce the concept of $H^s$-moderate, $H^s$-negligible nets, and a very weak solution of Sobolev order $s$. 
	\begin{definition}\label{H^k-moderate}
		\leavevmode
		
		\begin{itemize}
			\item[(i)] A net $(v_\eps)_\eps \in H^s(\mathbb{R}^n)^{(0,1]}$ is $H^s$-moderate if there exists $N \in \mathbb{N}_0$ and $c>0$ such that 
			\begin{equation*}
			|| v_\eps(x)||_{H^s(\mathbb{R}^n)} \leq c \eps^{-N},
			\end{equation*} 
			uniformly in $\eps \in (0,1]$.
			\item[(ii)] A net $(v_\eps)_\eps \in H^s (\mathbb{R}^n)^{(0,1]}$ is $H^s$-negligible if for all $q\in \mathbb{N}_0$ there exists $c>0$ such that 
			\begin{equation*}
			|| v_\eps(x)||_{H^s (\mathbb{R}^n)} \leq c\eps^q,  
			\end{equation*} 
			uniformly in $\eps \in (0,1]$.
		\end{itemize}
	\end{definition}
	Analogously, one can replace $H^s(\mathbb R^n)$ with $C([0,T], H^s(\mathbb R^n))$ and state the corresponding notions of $C([0,T], H^s(\mathbb R^n))$-moderate and $C([0,T], H^s(\mathbb R^n))$-negligible net. This will appear in the following notion of \emph{very weak solution of Sobolev order s}.
	\begin{definition}
		\label{def_vws_2}
		We say that the Cauchy problem 
		\[
		\begin{cases}
		D_t u=A(t,x,D_x)u+L(t,x,D_x)u+f(t,x),\\
		u(0,x) = g_{0} (x)
		\end{cases}
		\]
		admits a \emph{very weak solution} $(u_\eps)_\eps$ of \emph{Sobolev order s} if the net $(u_\eps)_\eps$ solves the regularised problem \eqref{Main-eq} and it is $C([0,T], H^s(\mathbb R^n))$-moderate.
	\end{definition}

	We are now ready to investigate the Cauchy problem \eqref{Main-eq} more in detail.


	\subsection{Hyperbolic systems with non-diagonalisable principal part}
	Consider the regularised hyperbolic system of type 
	\[
	\begin{cases}
	D_t u_{\varepsilon} = A_{\varepsilon} (t,x,D_x) u_{\varepsilon} + L_{\varepsilon} (t,x,D_x)u_{\eps} + f_{\varepsilon}(t,x), & \quad (t,x) \in [0,T]\times \mathbb{R}^n,\\
	u_{\varepsilon}(0,x) = g_{\varepsilon}(x), & \quad x \in \mathbb{R}^n,
	\end{cases}
	\]
	where $u_\eps(0,x) = [g_{1,\eps} (x), g_{2, \eps}(x)]^T$ and $f_\eps(t,x) = [f_{1,\eps}(t,x), f_{2, \eps}(t,x)]^T$, with the  $A_{\eps}(t,x,D_x)$ and $L_{\eps}(t,x,D_x)$ given by
	\begin{equation*}
	A_{\varepsilon}(t,x,D_x) = \begin{bmatrix}
	\sum_{j=1}^n 	\lambda_{1j, \varepsilon}(t,x) D_{x_j} & a_{12, \varepsilon}(t,x,D_x) \\
	0 & \sum_{j=1}^n \lambda_{2j, \varepsilon}(t,x)D_{x_j}
	\end{bmatrix},
	\end{equation*}
	and 
	\begin{equation*}
	L_{\varepsilon}(t,x,D_x) = \begin{bmatrix}
	\ell_{11, \varepsilon}(t,x) & \ell_{12, \varepsilon}(t,x) \\
	\ell_{21, \varepsilon}(t,x)\langle D_x\rangle^{-1} & \ell_{22, \varepsilon}(t,x)
	\end{bmatrix},
	\end{equation*}
	where $(\lambda_{1j, \eps})_\eps, (\lambda_{2j, \eps})_\eps \in  C([0,T], C^\infty(\mathbb{R}^n))^{(0,1]}$ for $j=1,\ldots, n$ and $(a_{12, \eps})_\eps$ and $(\ell_{ik, \eps})_\eps$ are nets of symbols of order $1$ and $0$ respectively, i.e., $(a_{12, \eps})_\eps \in C([0,T], S^1_{0,1}(\mathbb{R}^n \times \mathbb{R}^n ))^{(0,1]}$ and $(\ell_{ik, \eps})_\eps \in C([0,T], S^0_{0,1}(\mathbb{R}^n \times \mathbb{R}^n ))^{(0,1]}$ for $i,k=1,2$.
	
	The net of solutions $(u_\eps)_\eps$ of this hyperbolic system can be formulated in terms of operators $G_{i, \eps}^0$ and $G_{i, \eps}$ introduced in Subsection \ref{Aux_rem} if the following hypotheses hold:
	\begin{itemize}
		\item[(H1)] The coefficients $(\lambda_{1j, \eps})_\eps$ and $(\lambda_{2j, \eps})_\eps$ are real valued and $C([0,T], C^\infty(\mathbb{R}^n))$-moderate for $j=1,\ldots,n$ with compact support in $x$.
		\item[(H2)]  Also, $(\partial_k \lambda_{1j, \eps})_\eps$, $(\partial_k\lambda_{2j, \eps})_\eps$ and $(\ell_{ii, \eps})_{\eps}$ are logarithmic moderate for $k,j=1,\ldots,n$ and $i=1,2$.
	\end{itemize}
	Note that these assumptions allow to write $(u_{1,\eps})_\eps$ and $(u_{2,\eps})_\eps$ as in \eqref{eq:u21} and \eqref{eq:u22}. More details can be found in Theorem \ref{theo_FIO_inh}.

	
	In order to prove the existence of a very weak solution we start from the component $u_{1,\eps}$ and the equation
	\[
	u_{1,\eps} = \widetilde{U}^0_{1,\eps} + G_{1,\eps}(a_{12,\eps}G_{2,\eps}(\ell_{21,\eps}\langle D_x\rangle^{-1}u_{1,\eps})) + G_{1,\eps}(\ell_{12,\eps}G_{2,\eps}(\ell_{21,\eps}\langle D_x\rangle^{-1} u_{1,\eps})),
	\]
	where 
	\[
	\widetilde{U}^0_{1,\eps} = G_{1,\eps}^0 g_{1,\eps} + G_{1,\eps}(f_{1,\eps}) + G_{1,\eps}((a_{12,\eps}+\ell_{12,\eps})U^0_{2,\eps}).
	\]
	We make use of the Banach fixed point theorem assuming the additional hypothesis:
	\begin{itemize}
		\item[(H3)]  The operators $G^0_{i,\eps}$ and $G_{i,\eps}$ map $H^s$-moderate and $C([0,T], H^s(\mathbb{R}^n))$-moderate nets into themselves, respectively and 
		\[\mathcal{G}^0_{1, \eps} := G_{1, \eps} \circ (a_{12, \eps} + \ell_{12, \eps})\circ G_{2, \eps} \circ \ell_{21, \eps}\langle D_x\rangle^{-1}
		\]
		has the operator norm in $H^s$ strictly less than 1
		\begin{equation*}
		|| \mathcal{G}^0_{1,\eps}  ||_{H^s \rightarrow H^s} = O(1).
		\end{equation*}
		Note that $\mathcal{G}^0_{1,\eps} $ acts continuously on $H^s$ since it is of order $0$. 
	\end{itemize} 
	Let us now work on better understanding the hypothesis (H3). We have that (H3) holds if    
	\begin{align*}
	|| G_{1, \eps} \circ a_{12, \varepsilon} \circ G_{2, \eps} \circ \ell_{21, \varepsilon} \langle D_x\rangle^{-1}||_{H^s \rightarrow H^s} \leq T C_1(\varepsilon) A_1(\varepsilon) a_{12}(\varepsilon) C_2(\varepsilon) A_2(\eps) L_{21}(\varepsilon)=O(1) , 
	\end{align*}
	and 
	\begin{align*}
	|| G_{1, \eps} \circ \ell_{12, \varepsilon} \circ G_{2, \eps} \circ \ell_{21, \varepsilon}\langle D_x\rangle^{-1}||_{H^s \rightarrow H^s } \leq T C_1(\varepsilon) A_1(\varepsilon) L_{12}(\varepsilon) C_2(\varepsilon) A_2(\eps) L_{21}(\varepsilon)=O(1) , 
	\end{align*}
	where 
	\begin{align*}
	A_i(\varepsilon) &:= \sup_{t,s\in[0,T]}\sup_{|\beta|\leq 2n+1} || \partial_x^{\beta} A_{i, \varepsilon}(s,t,x)||_{L^{\infty}(\mathbb{R}^n_x)} \,\, \text{ for }\,\, i=1,2, \\
	a_{12}(\varepsilon) &:=  \sup_{t\in[0,T]}\sup_{|\alpha|,|\beta|\leq 2n+1} || \partial_\xi^{\alpha} \partial_{x}^{\beta}a_{12, \varepsilon}(t,x,\xi)||_{L^{\infty}(\mathbb{R}^n_x\times \mathbb{R}^n_{\xi})}, \\
	L_{12}(\varepsilon) &:= \sup_{t\in[0,T]}\sup_{|\beta|\leq 2n+1} || \partial_x^{\beta} \ell_{12, \varepsilon}(t,x)||_{L^{\infty}(\mathbb{R}^n_x)},\\
	L_{21}(\varepsilon) &:= \sup_{t\in[0,T]}\sup_{|\beta|\leq 2n+1} || \partial_{x}^{\beta}\ell_{21, \varepsilon}(t,x)||_{L^{\infty}(\mathbb{R}^n_x)},
	\end{align*}
	and $C_1(\varepsilon)$ and $C_2(\varepsilon)$ include the combination of constants $C(\eps)$, $C_{\alpha, i}(\eps)$ and $ C_{\beta, i}(\eps)$ from
	\begin{align*}
	|\partial_{\xi} \phi_{i, \eps}(t, x,\xi)-\partial_{\xi}\phi_{i, \eps}(t, y,\xi)| &\geq C(\eps)|x-y| \,\,\,\, \text{ for } \,\, x,y\in \mathbb{R}^n, \,\, \xi \in \mathbb{R}^n,  \\
	|\partial_y\phi_{i, \eps}(t, y,\xi)-\partial_y \phi_{i, \eps}(t, y,\eta)| &\geq C(\eps) |\xi - \eta| \,\,\,\, \text{ for } \,\,  y\in \mathbb{R}^n, \,\, \xi, \eta \in \mathbb{R}^n, \nonumber
	\end{align*}
	and that 
	\begin{equation*}
	|\partial_y^{\alpha}\partial_{\xi}\phi_{i, \eps}(t, y,\xi)| \leq C_{\alpha, i}(\eps), \,\, 	|\partial_y\partial_{\xi}^{\beta}\phi_{i, \eps}(t, y,\xi)| \leq C_{\beta, i}(\eps),
	\end{equation*}
	for $1\leq |\alpha|, |\beta|\leq 2n+2$ and $i=1,2$.
	
	We refer the reader to the appendix at the end of the paper for the Sobolev mapping properties of Fourier integral operators that we have employed above. Note that the nets of constants $C_{\alpha, i}(\eps)$, $C_{\beta, i}(\eps)$ are automatically moderate while $C(\eps)$ and $C^{-1}(\eps)$ are assumed to be moderate to make sure that our operators have the right Sobolev mapping properties. The stronger hypothesis (H3) is required to allow a fixed point argument independent of the parameter $\varepsilon$. We therefore conclude that (H3) can be written as
	
	\begin{itemize}
		\item[(H3)] the nets of constants $(C(\eps))_\eps$ and $(C^{-1}(\eps))_\eps$ are moderate and \begin{equation*}
		|| \mathcal{G}^0_{1,\eps}  ||_{H^s \rightarrow H^s} = O(1).
		\end{equation*}
	\end{itemize}
	
	If (H3) holds, then we will apply Banach's fixed point theorem in the space $X(t):= C([0,t], H^s(\mathbb{R}^n))$ for $t \in [0,T]$ with the norm 
	\begin{equation*}
	|| u_{1, \eps}||_{X(t)} = \sup_{0\leq \tau \leq t} || u_{1,\eps}(\tau, \cdot)||_{H^s}. 
	\end{equation*} 
	Note that we can rewrite the equation for $u_{1,\eps}$ as 
	\begin{equation*}
	u_{1,\eps} = {\tilde{U}}_{1, \eps}^0 + \mathcal{G}^0_{1,\eps} u_{1,\eps}.
	\end{equation*}
	
	Using composition of Fourier integral operators and hypothesis (H3) we have that the 0-order Fourier integral operator $\mathcal{G}^0_{1, \eps}$ maps $X(t)$ continuously into itself and for small time interval it is a contraction, in the sense that there exists $T^*\in [0,T]$ such that 
	\begin{equation}
	|| \mathcal{G}^0_{1, \eps} (u-v)||_{X(T^*)} \leq C T^* || u-v||_{X(T^*)},
	\end{equation}
	with $C T^*<1$. The existence of a unique fixed point $u_{1,\eps}$ for the map $\mathcal{G}^0_{1,\eps}$ is provided by Banach's fixed point theorem and is equivalent to say that the operator $I-\mathcal{G}^0_{1, \eps}$ is invertible on $[0, T^*]$ for all values of $\eps$. It therefore follows that $(u_{1,\eps})_\eps$ inherits the moderateness properties of $({\tilde{U}}_{1, \eps}^0)_\eps$. Note that 
	as already observed in \cite{GarJRuz} the constant $C$ does not depend on the initial data so the argument can be iterated to cover the full interval $[0, T]$.
	
	We have therefore proven the following general result.

	\begin{theorem}\label{thm_general_statmenet}
		Let us consider the Cauchy problem 
		\begin{equation}\label{eq-thm}
		\begin{cases}
		D_t u = A (t,x, D_x)u + L (t,x,D_x) u + f(t,x), & \quad (t,x) \in [0,T]\times \mathbb{R}^n,\\
		u(0,x)= g(x), & \quad x \in \mathbb{R}^n,
		\end{cases}
		\end{equation}
		where  
		\begin{equation*}
		A(t,x,D_x) = \begin{bmatrix}
		\sum_{j=1}^n 	\lambda_{1j}(t,x) D_{x_j} & a_{12}(t,x,D_x) \\
		0 & \sum_{j=1}^n \lambda_{2j}(t,x)D_{x_j}
		\end{bmatrix}
		\end{equation*}
		is an upper triangular matrix of first order differential operator and 
		\begin{equation*}
		L(t,x,D_x) = \begin{bmatrix}
		\ell_{11}(t,x) & \ell_{12}(t,x) \\
		\ell_{21}(t,x)\langle D_x\rangle^{-1} & \ell_{22}(t,x)
		\end{bmatrix},
		\end{equation*}
		is a matrix of pseudo-differential operator of order 0 with $\ell_{12}$ of order $-1$, continuous with respect to $t$ and such that 
		\begin{itemize}
			\item[(i)] all coefficients of $A(t,x,D_x)$ and $L(t,x,D_x)$ are in $C([0,T], \mathcal{E}'(\mathbb{R}^n))$ with compact support in $x$,
			\item[(ii)] the initial data $g_i(x)\in \mathcal{E}'(\mathbb{R}^n)$ and the source term $f_i(t,x) \in C([0,T], \mathcal{E}'(\mathbb{R}^n))$ with compact support in $x$ for $i=1,2$. 
		\end{itemize}
		
		Assume that (H1), (H2), (H3) are satisfied. Let the nets of regularised initial data $(g_{i, \eps})_\eps$ and the right-hand side $(f_{i, \eps})_\eps$ be $H^{s +i -1}$-moderate and $C([0,T], H^{s+i-1})$-moderate for $i=1,2$, respectively. Then the net of solutions $(u_{\eps})_\eps$ is a very weak solution of anisotropic Sobolev type, i.e. $u_{i,\eps}$ is $C([0,T], H^{s+i-1})$-moderate for $i=1,2$.
	\end{theorem}
	
	\begin{remark}
		Note that when regularising a distribution with compact support with a mollifier which is also compactly supported we automatically get a moderate net of smooth functions with compact support and therefore a moderate net of any Sobolev order.
	\end{remark}

	The existence of a very weak solution in the sense of Definition \ref{def_vws_1} follows from the following corollary.
	
	\begin{corollary}
		\label{cor_vws}
		If (H1) and (H2), (H3) hold then the Cauchy problem \eqref{eq-thm} has a very weak solution $(u_\eps)_\eps$, i.e., a net $(u_\eps)_\eps$ which solves the regularised Cauchy problem and it is $C([0,T], C^\infty(\mathbb{R}^n))$-moderate.  
	\end{corollary}
	
	\begin{remark}
		\label{rem_uniqueness}
		Negligible perturbations of the system coefficients and the initial data leads to a negligible perturbation in the solution. This follows from the fact that our nets of operators maps negligible nets into negligible nets. For this reason we say that our Cauchy problem is well-posed in a very weak sense.
	\end{remark}
	
	\subsection{Special cases of the hyperbolic system in \eqref{eq-thm}}
	In this subsection, we analyse some special cases of Theorem \ref{thm_general_statmenet}.
	First, we assume that the matrix $L(t,x,D_x)$ of the lower order terms is  upper triangular and that $A(t,x,D_x)$ is either diagonal or in upper triangular form. In both these cases the hypotheses (H1) and (H2) are sufficient to prove Theorem \ref{thm_general_statmenet}.
	\begin{corollary}\label{cor-1}
		Let us consider the Cauchy problem \eqref{eq-thm} with a diagonal matrix of first order differential operators $A(t,x,D_x)$ and an upper triangular matrix of 0-order pseudo-differential operator  $L(t,x,D_x)$ as
		\begin{equation*}
		A(t,x,D_x) = \begin{bmatrix}
		\sum_{j=1}^n 	\lambda_{1j}(t,x) D_{x_j} & 0 \\
		0 & \sum_{j=1}^n \lambda_{2j}(t,x)D_{x_j}
		\end{bmatrix},
		\end{equation*}
		and 
		\begin{equation*}
		L(t,x,D_x) = \begin{bmatrix}
		\ell_{11}(t,x) & \ell_{12}(t,x) \\
		0 & \ell_{22}(t,x)
		\end{bmatrix},
		\end{equation*}
		where
		\begin{itemize}
			\item[(i)] all coefficients of $A(t,x,D_x)$ and $L(t,x,D_x)$ are in $C([0,T], \mathcal{E}'(\mathbb{R}^n))$ with compact support in $x$,
			\item[(ii)] the initial data $g_i(x)\in \mathcal{E}'(\mathbb{R}^n)$ and the source term $f_i(t,x) \in C([0,T], \mathcal{E}'(\mathbb{R}^n))$ with compact support in $x$ for $i=1,2$. 
		\end{itemize}
		
		Assume that (H1) and (H2) are satisfied. 
		Then the net of solutions $(u_{\eps})_\eps$ of the regularised problem is a very weak solution, i.e. $(u_{i,\eps})_\eps$ is $C([0,T], C^\infty(\mathbb{R}^n))$-moderate for $i=1,2$ and can be written as
		\begin{align*}
		u_{1, \varepsilon} &= G_{1,\eps}^0 g_{1, \varepsilon} + G_{1, \eps} f_{ 1, \varepsilon} ,\\
		u_{2, \varepsilon} &= G_{2, \eps}^0 g_{2, \varepsilon} + G_{2, \eps} f_{2, \varepsilon}.
		\end{align*}
	\end{corollary}
	\begin{corollary}\label{cor-2}
		Let the matrices $A(t,x,D_x)$ and $L(t,x,D_x)$ in the Cauchy problem \eqref{eq-thm} be both upper triangular.  
		Let
		\begin{itemize}
			\item[(i)] all coefficients of $A(t,x,D_x)$ and $L(t,x,D_x)$ are in $C([0,T], \mathcal{E}'(\mathbb{R}^n))$ with compact support in $x$,
			\item[(ii)] the initial data $g_i(x)\in \mathcal{E}'(\mathbb{R}^n)$ and the source term $f_i(t,x) \in C([0,T], \mathcal{E}'(\mathbb{R}^n))$ with compact support in $x$ for $i=1,2$. 
		\end{itemize}
		
		Assume that (H1) and (H2) are satisfied. Assume that there exist moderate nets of constants $(C(\eps))_\eps$ and $(C^{-1}(\eps))_\eps$ such that 
		\begin{align}
		\label{est-C-eps}
		|\partial_{\xi} \phi_{i, \eps}(t, x,\xi)-\partial_{\xi}\phi_{i, \eps}(t, y,\xi)| &\geq C(\eps)|x-y| \,\,\,\, \text{ for } \,\, x,y\in \mathbb{R}^n, \,\, \xi \in \mathbb{R}^n,  \\
		|\partial_y\phi_{i, \eps}(t, y,\xi)-\partial_y \phi_{i, \eps}(t, y,\eta)| &\geq C(\eps) |\xi - \eta| \,\,\,\, \text{ for } \,\,  y\in \mathbb{R}^n, \,\, \xi, \eta \in \mathbb{R}^n, \nonumber
		\end{align}
		for all $i=1,2$. Then the net of solutions $(u_{\eps})_\eps$ of the regularised problem is a very weak solution of anisotropic Sobolev type for all $s\in\mathbb{R}$, i.e. $(u_{i,\eps})_\eps$ is $C([0,T], H^{s+i-1})$-moderate for $i=1,2$ and can be written as
		\begin{align*}
		u_{1, \varepsilon} &=  G_{1,\eps}^0 g_{1, \varepsilon} + G_{1, \eps} f_{ 1, \varepsilon}  +G_{1, \eps} (a_{12,\eps} U^0_{2,\eps}),\\
		u_{2, \varepsilon} &=U^0_{2,\eps} = G_{2, \eps}^0 g_{2, \varepsilon} + G_{2, \eps} f_{2, \varepsilon},
		\end{align*}
		where $U^0_{2,\eps}$ is defined in \eqref{eq:aux52}.
	\end{corollary}
	\begin{remark}
		It is important to note that Corollary \ref{cor-1} is independent of the hypothesis (H3) and that Corollary \ref{cor-2} only requires the hypothesis on the phase function which guarantees the right Sobolev mapping properties for our Fourier integral operators. This is due to the fact that the Banach fixed point theorem is not required for the proof of Corollaries \ref{cor-1} and \ref{cor-2}. It is also clear that Corollary \ref{cor-1} and \ref{cor-2} hold for $m\times m$-systems as well. It is indeed sufficient to formulate the hypotheses (H1) and (H2) with $i=1,\ldots, m$ and assume that \eqref{est-C-eps} holds for $i=1,\ldots, m$.
	\end{remark}
	
	We now focus on $2\times 2$-systems where the hypothesis (H3) is needed to prove the existence of a very weak solution. For the sake of the reader we reformulate Proposition \ref{Prop-L2} for the operators involved in (H3). This will apply to a scenario where the matrices $A(t,x,D_x)$ and $L(t,x,D_x)$ are not as in the previous two corollaries. Let $\phi_{i,\eps}(t,x,\xi)=\gamma_{i,\eps}(x,t;0)\xi$. We have that if there exists constants $C, C_\alpha, C_\beta$ independent of $\eps\in(0,1]$ such that 
	
	\begin{align*}
	| \gamma_{i,\eps}(t, x; 0) -\gamma_{i,\eps}(t, y; 0)| &\geq C|x-y| \,\,\,\, \text{ for } \,\, t\in[0,T],\, x,y\in \mathbb{R}^n, \\
	|\partial_y\gamma_{i, \eps}(t, y; 0)\cdot (\xi-\eta)| &\geq C |\xi - \eta| \,\,\,\, \text{ for } \,\,  t\in[0,T],\, y\in \mathbb{R}^n, \,\, \xi, \eta \in \mathbb{R}^n, \nonumber
	\end{align*}
	and 
	\begin{equation*}
	|\partial_y^{\alpha}\gamma_{i, \eps}(t, y; 0)| \leq C_{\alpha}, \,\, 	|\partial_y\partial_{\xi}^{\beta}(\gamma_{i,\eps}(t, y; 0)\cdot \xi)| \leq C_{\beta},\,\,\, \text{for}\, t\in[0,T],\, y\in\mathbb{R}^n,
	\end{equation*}
	for $1\leq |\alpha|, |\beta|\leq 2n+2$ and $i=1,2$,
	and the lower order terms $\ell_{12,\eps}$ and $\ell_{21,\eps}$ have derivatives up to order $2n+1$ which are bounded with respect to the parameter $\eps$ as well, then the hypothesis (H3) is automatically fulfilled. This corresponds to a case of higher regularity of the coefficients which is treated in the next corollary.

	\begin{corollary}\label{cor-3}
		Let us consider the Cauchy problem \eqref{eq-thm} with given a diagonal matrix of first order differential operators $A(t,x,D_x)$ and a matrix of 0-order pseudo-differential operator  $L(t,x,D_x)$ as
		\begin{equation*}
		A(t,x,D_x) = \begin{bmatrix}
		\sum_{j=1}^n 	\lambda_{1j}(t,x) D_{x_j} & 0 \\
		0 & \sum_{j=1}^n \lambda_{2j}(t,x)D_{x_j}
		\end{bmatrix},
		\end{equation*}
		and 
		\begin{equation*}
		L(t,x,D_x) = \begin{bmatrix}
		0 & \ell_{12}(t,x) \\
		\ell_{21}(t,x)  & 0
		\end{bmatrix},
		\end{equation*}
		where
		\begin{itemize}
			\item[(i)] the diagonal elements $\lambda_{1j}$, $\lambda_{2j} \in C([0,T], C^{k+1}(\mathbb{R}^n))$ for $j=1,\ldots, n$ and $k\geq 2n+1$ with compact support in $x$;
			\item[(ii)] the initial data $g(x)\in \mathcal{E}'(\mathbb{R}^n)$ and the source term $f_i(t,x) \in C([0,T], \mathcal{E}'(\mathbb{R}^n))$ with compact support in $x$, for $i=1,2$;
			\item[(iii)] the non-diagonal elements $\ell_{12}(t,x), \ell_{21}(t,x) \in C([0,T], C^k(\mathbb R^n))$ with compact support with respect to $x$. 
		\end{itemize}
		
		Then the net of solutions $(u_{\eps})_\eps$ of the regularised problem is a very weak solution of Sobolev type for all $s\in\mathbb{R}$, i.e. $(u_{i,\eps})_\eps$ is $C([0,T], H^{s})$-moderate for $i=1,2$.
	\end{corollary}
	
	\begin{proof}[Proof of Corollary \ref{cor-3}]
		The enhanced regularity of the coefficients allows (H1), (H2), and (H3) to be satisfied automatically with constants that do not depend on $\eps$. Under assumptions (H1) and (H2), we can express the $(u_{\eps})_{\eps}$ components as:
		\begin{align*}
		u_{1, \eps} &= G^0_{1, \eps} g_{1, \eps} + G_{1, \eps} f_{1, \eps} + G_{1,\eps}(\ell_{12, \eps}u_{2, \eps}), \\
		u_{2, \eps}& = G^0_{2, \eps} g_{2, \eps} + G_{2, \eps} f_{2, \eps} + G_{2,\eps}(\ell_{21, \eps}u_{1, \eps}).
		\end{align*}
		By substituting  $u_{2, \eps}$ into the equation of $u_{1,\eps}$, we get
		\begin{equation*}
		u_{1,\eps} = \widetilde{U}_{1,\eps}^0 + \mathcal{G}^0_{1,\eps} u_{1,\eps},
		\end{equation*}
		where
		\begin{align*}
		\widetilde{U}_{1,\eps}^0 &:=  G^0_{1, \eps} g_{1, \eps} + G_{1, \eps} f_{1, \eps}+ G_{1,\eps}(\ell_{12} \circ G^0_{2,\eps}g_{2, \eps}) + G_{1,\eps}(\ell_{12} \circ G_{2,\eps}f_{2, \eps}), \\
		\mathcal{G}^0_{1,\eps} &:= G_{1,\eps} \circ \ell_{12, \eps} \circ G_{2,\eps}\circ \ell_{21, \eps}. 
		\end{align*}
		Note that $\mathcal{G}^0_{1, \eps}$ is a zero-order operator. In other words it is not necessary in this specific case to assume that the entry below the diagonal in the matrix $L(t,x,D_x)$ is of order $-1$. Thus, one can prove that $\mathcal{G}^0_{1, \eps}$ has the operator norm in $H^s$ strictly less than 1 due to (H3). Moreover, $\widetilde{U}^0_{1, \eps} $ is $H^s$-moderate. By applying the Banach fixed-point theorem, we prove the existence of very weak solution $(u_{1,\eps})_\eps$, i.e., $(u_{1,\eps})_\eps$ is $C([0,T], H^s)$-moderate. Finally, we obtain $u_{2, \eps}$ by substituting $u_{1, \eps}$ into the equation for $u_{2, \eps}$.
	\end{proof}
	
	\section{Consistency and Applications}
	In this section, we prove that our result is consistent with the classical well-posedness result obtained in \cite{GarJRuz} when the system is regular enough. In detail, we show that every very weak solution converges to the classical solution when it e\-xi\-sts, namely when the system has smooth coefficients with respect to $x$. In addition, we discuss some examples of physical relevance. 
	
	\begin{theorem}\label{thm_cnstnc}
		Consider the Cauchy problem \eqref{eq-thm}, where $A(t,x,D_x)$ is an upper triangular matrix of first-order pseudo-differential operators and $L(t,x,D_x)$ is a matrix of zero-order pseudo-differential operators, continuous with respect to $t$, of the form
		\begin{equation*}
		A(t,x,D_x) = \begin{bmatrix}
		\sum_{j=1}^n 	\lambda_{1j}(t,x) D_{x_j} & a_{12}(t,x,D_x) \\
		0 & \sum_{j=1}^n \lambda_{2j}(t,x)D_{x_j}
		\end{bmatrix}
		\end{equation*}
		and
		\begin{equation*}
		L(t,x,D_x) = \begin{bmatrix}
		\ell_{11}(t,x) & \ell_{12}(t,x) \\
		\ell_{21}(t,x)\langle D_x\rangle^{-1} & \ell_{22}(t,x)
		\end{bmatrix}.
		\end{equation*}
		Assume that all coefficients of $A(t,x,D_x)$ and $L(t,x,D_x)$ are continuous in $t$, smooth and compactly supported with respect to $x$, the initial data $g_i(x) \in C^\infty_c(\mathbb{R}^n)$ and $f_i(t,x) \in C([0,T], C^\infty_c(\mathbb{R}^n))$ for $i=1,2$ and $s\in\mathbb{R}$. Hence, 
		\begin{itemize}
			\item[(i)] the Cauchy problem \eqref{eq-thm}  has a unique solution $u$ with components $u_i\in C([0,T], C^\infty_c(\mathbb{R}^n))$ for $i=1,2$;
			\item[(ii)]  the Cauchy problem \eqref{eq-thm} has a very weak solution $(u_{\eps})_{\eps}$, i.e., the components $(u_{i,\eps})_\eps$ are $C([0,T], C^\infty(\mathbb{R}^n))$-moderate for $i=1,2$.
			\item[(iii)] the net $u_\eps(t,\cdot)$ converges to $u(t,\cdot)$ as $\eps\to 0$ in $L^2(\mathbb{R}^n_x)$ uniformly with respect to $t \in [0, T ]$.
		\end{itemize}
	\end{theorem}
	
	\begin{proof}[Proof of Theorem \ref{thm_cnstnc}]
		(i) The well-posedness of the Cauchy problem can be derived from the results presented by the first author, J\"ah, and Ruzhansky in \cite{GarJRuz}.
		
		(ii) The existence of a very weak solution of anisotropic Sobolev type is obtained from Theorem \ref{thm_general_statmenet}. Note that since the coefficients are smooth the hypotheses (H1)-(H3) are fulfilled with constants independent of $\eps$.
		
		(iii) Let $(u_\eps)_\eps$ and $u$ be a very weak solution and the classic solution of our Cauchy problem, respectively. Our argument is independent of the choice of the mollifier and the scale we are using in our regularisation, so we namely prove that every very weak solution produced in this way is convergent to the unique classical solution.
		
		We now compare the classical Cauchy problem with solution $u$, with the regularised one with solution $(u_\eps)_\eps$. It follows that
		\begin{equation*}
		\begin{cases}
		D_t (u_\eps -u) = A_\eps (u_\eps-u) + L_\eps (u_\eps - u) + (f_\eps - f) + R_\eps, & \quad (t,x) \in [0,T]\times \mathbb{R}^n,\\
		(u_\eps - u)(0,x)= (g_\eps - g)(x), & \quad x \in \mathbb{R}^n,
		\end{cases}
		\end{equation*}
		where 
		\begin{equation*}
		R_\eps(t,x) = [(A_\eps - A)(t,x,D_x) + (L_\eps - L)(t,x,D_x)] u(t,x).
		\end{equation*}
		Repeating the argument in Theorem \ref{thm_general_statmenet},  we arrive at
		\begin{align*}
		V_{1,\eps} = U^0_{1,\eps} + G_{1, \eps}\circ (a_{12, \eps} + \ell_{12, \eps}) \circ G_{2,\eps} \circ \ell_{21, \eps} \langle D_x \rangle^{-1} V_{1,\eps},    
		\end{align*}
		where 
		\begin{align*}
		V_{1,\eps} &= u_{1,\eps} - u_1, \\
		U^0_{1, \eps} &= G_{1, \eps}^0 (g_{1, \eps} - g_1) + G_{1,\eps} ([f_{1,\eps} - f_1] + R_{1,\eps}) \\
		&+ G_{1, \eps} \circ (a_{12, \eps} + \ell_{12, \eps}) \circ [G_{2,\eps}^0 (g_{2, \eps}-g_2) + G_{2,\eps}(f_{2, \eps}-f_2) + G_{2, \eps} R_{2,\eps}],
		\end{align*}
		and 
		\begin{align*}
		R_{1,\eps} &:= [(\lambda_{1,\eps} - \lambda_1)(t,x,D_x)u_1 +(a_{12, \eps} - a_{12})u_2 + (\ell_{11, \eps} - \ell_{11})u_1 + (\ell_{12, \eps} - \ell_{12})u_2], \\
		R_{2,\eps} &:= [(\lambda_{2,\eps} - \lambda_2)(t,x,D_x)u_2 + (\ell_{21, \eps} - \ell_{21})\langle D_x\rangle^{-1} u_1 + (\ell_{22, \eps} - \ell_{22})u_2].
		\end{align*}
		By using the hypothesis (H3) and the Banach fixed-point theorem with $CT \leq 1$, we get 
		\begin{align*}
		||V_{1,\eps}||_{L^2(\mathbb{R}^n_x)} &= ||U^0_{1,\eps}||_{L^2(\mathbb{R}^n_x)} + ||G_{1, \eps}\circ (a_{12, \eps} + \ell_{12, \eps}) \circ G_{2,\eps} \circ \ell_{21, \eps}\langle D_x\rangle^{-1}  V_{1,\eps}||_{L^2(\mathbb{R}^n_x)} \\  
		& \leq ||U^0_{1,\eps}||_{L^2(\mathbb{R}^n_x)} + CT ||V_{1,\eps}||_{L^2(\mathbb{R}^n_x)} \\
		& \leq ||U^0_{1,\eps}||_{L^2(\mathbb{R}^n_x)}. 
		\end{align*}
		Note that because of the regularity of the coefficients the operators involved above fulfill estimates which are independent of the parameter $\eps\in (0,1]$. We therefore have that $||U^0_{1,\eps}||_{L^2(\mathbb{R}^n_x)} \rightarrow 0$ as $\eps \rightarrow 0$ since $|| g_{1,\eps} - g_1||_{L^2(\mathbb{R}^n_x)} \rightarrow 0$, $|| f_{1, \eps} - f_1 ||_{L^2(\mathbb{R}^n_x)} \rightarrow 0$, $|| \lambda_{i, \eps} - \lambda_i||_{L^\infty(\mathbb{R}^n_x)}\rightarrow 0$, $\langle{x}\rangle^{-1}|| a_{12, \eps} - a_{12}||_{L^\infty(\mathbb{R}^n_x)}\rightarrow 0$, and $|| \ell_{ij, \eps} - \ell_{ij} ||_{L^\infty(\mathbb{R}^n_x)}\rightarrow 0$ for $i,j =1,2$ as $\eps \rightarrow 0$. This gives that $|| u_{1,\eps} - u_1 ||_{L^2(\mathbb{R}^n_x)} \rightarrow 0$. From the continuity assumption we have that all the limits above hold uniformly with respect to $t\in[0,T]$. This yields to $\Vert u_\eps-u\Vert_{L^2(\mathbb{R}^n_x)}\to 0$ as $\eps\to 0$ uniformly with respect to $t\in[0,T]$.
	\end{proof}
	We conclude this section with the analysis of few examples where the system coefficients are less than continuous or in general distributions with compact support.

	\subsection{Examples}
	For the sake of simplicity we work in space dimension 1. Let us consider the Cauchy problem 
	\begin{equation*}
	\begin{cases}
	D_t u =\begin{bmatrix}
	H(x)D_x  & a_{12}(t,x,D_x) \\
	0 & H(x)D_x
	\end{bmatrix} u + \begin{bmatrix}
	\ell_{11}(t,x) & \ell_{12}(t,x) \\
	\ell_{21}(t,x)\langle D_x\rangle^{-1}& \ell_{22}(t,x)
	\end{bmatrix} u,\\
	u(0,x) = g(x),  
	\end{cases}
	\end{equation*}
	where $(t,x) \in [0,T]\times \mathbb{R}$, $u(0,x) = [g_{1} (x), g_{2}(x)]^T$ with $g_1, g_2 \in \mathcal{E}'(\mathbb R)$ and the operators $a_{12}$ and $\ell_{ij}$ are pseudo-differential operators of order $1$ and $0$ respectively continuous with respect to $t$.
	For the moment we assume that only the coefficients on the diagonal of the matrix $A$ are non-regular, i.e., $H$ is the Heaviside function ($H(x) =1$ for $x\geq 0$ and $H(x)=0$ for $x<0$). One could in principle replace $H$ with any discontinuous but bounded function. Note that already in this situation the system is not treatable within the classical theory of hyperbolic systems with multiplicities that as in \cite{GarJRuz, GarJRuz2} requires smoothness with respect to the variable $x$. It is therefore meaningful to look at the regularised problem and at the net $(u_{\eps})_{\eps}$ of its solutions. In detail,    
	
	\begin{equation*}
	\begin{cases}
	D_t u_\eps =\begin{bmatrix}
	H_\eps(x)D_x  & a_{12}(t,x,D_x) \\
	0 & H_\eps(x)D_x
	\end{bmatrix} u_\eps + \begin{bmatrix}
	\ell_{11}(t,x) & \ell_{12}(t,x) \\
	\ell_{21}(t,x)\langle D_x\rangle^{-1}& \ell_{22}(t,x)
	\end{bmatrix} u_\eps,\\
	u_\eps(0,x) = g_\eps(x), 
	\end{cases}
	\end{equation*}
	where 
	\begin{align*}
	H_{\eps}(x) &= (H * \psi_{\omega(\eps)})(x),\\
	g_{1,\eps}(x) & = (g_1 * \psi_{\eps})(x), \\ 
	g_{2,\eps}(x) & = (g_2 * \psi_{\eps})(x).
	\end{align*}
	We will choose the scale $\omega(\eps)$ later in order to make sure that the hypothesis (H2) holds.
	
	We can reformulate the regularised problem in terms of Fourier integral operators follows
	\begin{align*}
	u_{1, \eps} &= G_{1, \eps}^0 g_{1, \eps} + G_{1, \eps} ((a_{12} + \ell_{12})u_{2, \eps} ),\\
	u_{2, \eps} &=  G_{2, \eps}^0 g_{2, \eps} + G_{2, \eps} ( \ell_{21}\langle D_x \rangle^{-1}u_{1, \eps} ).
	\end{align*}
	where 
	\begin{align*}
	G_{j, \eps}^0 g_{j, \eps}(t,x) &= \int_{\mathbb R} {\rm e}^{i (x + H_{\eps}(x)t)\xi } {\rm e}^{\int_0^t \ell_{jj}(\tau,x + H_{\eps}(x)\tau) d\tau } \hat{g}_{j, \eps}(\xi) d\xi\\
	&= {\rm e}^{\int_0^t \ell_{jj}(\tau,x + H_{\eps}(x)\tau) d\tau } g_{j, \eps}(x + H_{\eps}(x)t),\\
	G_{j, \eps} f_{ \eps}(t,x) &= \int_{0}^t \int_{\mathbb R} {\rm e}^{i (x + H_{\eps}(x)t)\xi } {\rm e}^{\int_s^t \ell_{jj}(\tau,x + H_{\eps}(x)\tau) d\tau } \hat{f}_{\eps}(s,\xi) d\xi ds\\
	&=\int_0^t\int_{\mathbb{R}^2}{\rm e}^{i (x + H_{\eps}(x)t)\xi-iy\xi}A_{j,\eps}(s,t,x)f_\eps(s,y)\, dy\dslash\xi ds\\
	&= \int_0^t {\rm e}^{\int_s^t \ell_{jj}(\tau,x + H_{\eps}(x)\tau) d\tau }  f_{ \eps}(s, x+ H_{\eps}(x)(t-s)) ds,
	\end{align*}
	for $j=1,2$.
	
	Let us now set $\omega^{-1}(\eps)=\log(\eps^{-1})$. Note that in this case the assumptions (H1) and (H2) are immediately fulfilled. Indeed, $(H_\eps)_\eps$ is real-valued and moderate and $H'_\eps=\psi_{\omega(\eps)}$ is logarithmic moderate. In addition $\Vert H_\eps\Vert_{L^\infty}\le 1$.  
	
	In the special case when both the matrices of the principal part and the lower order terms are diagonal the assumption (H3) is not needed. In detail, if $a_{12} =0$ and $\ell_{12}=\ell_{21}=0$, then   
	
	\begin{align*}
	u_{1, \eps} &= e^{\int_0^t \ell_{11}(\tau,x + H_{\eps}(x)\tau,D_x) d\tau } g_{1, \eps}(x + H_{\eps}(x)t), \\
	u_{2, \eps} &= e^{\int_0^t \ell_{22}(\tau,x + H_{\eps}(x)\tau,D_x) d\tau } g_{2, \eps}(x + H_{\eps}(x)t).
	\end{align*}

	In general for a full matrix $L$ one would need the assumption (H3) to be fulfilled as well. This is investigated in the following example where we also handle singularities in time. 
	
	Let
	\begin{equation*}
	\begin{cases}
	D_t u =\begin{bmatrix}
	\lambda_1(t)D_x  & a_{12}(t,x,D_x) \\
	0 & \lambda_2(t)D_x
	\end{bmatrix} u + \begin{bmatrix}
	\ell_{11}(t,x) & \ell_{12}(t,x) \\
	\ell_{21}(t,x)\langle D_x\rangle^{-1}& \ell_{22}(t,x)
	\end{bmatrix} u,\\
	u(0,x) = g(x),  
	\end{cases}
	\end{equation*}
	where $(t,x) \in [0,T]\times \mathbb{R}$, $\lambda_i$, $i=1,2$, are distributions with compact support in $[0,T]$, and the coefficients of the matrix of the lower order terms are distributions with compact support in $t$ and smooth in $x$. For the sake of simplicity we assume that $a_{12}\in C([0,T], S^1(\mathbb{R}^{2n}))$. As usual we have that $u(0,x) = [g_{1} (x), g_{2}(x)]^T$ with $g_1, g_2 \in \mathcal{E}'(\mathbb R)$. Instead of regularising in $x$ this time we regularise in $t$ and restrict the domain of our nets to the interval $[0,T]$ after regularisation. We obtain
	\begin{equation*}
	\begin{cases}
	D_t u_\eps =\begin{bmatrix}
	\lambda_{1,\eps}(t)D_x  & a_{12}(t,x,D_x) \\
	0 & \lambda_{2,\eps}(t)D_x
	\end{bmatrix} u_\eps + \begin{bmatrix}
	\ell_{11,\eps}(t,x) & \ell_{12,\eps}(t,x) \\
	\ell_{21,\eps}(t,x)\langle D_x\rangle^{-1}& \ell_{22,\eps}(t,x)
	\end{bmatrix} u_\eps,\\
	u_\eps(0,x) = g_\eps(x),  
	\end{cases}
	\end{equation*}
	where
	\[
	\lambda_{i,\eps}(t)=\lambda_i\ast\psi_{\eps},
	\]
	for $t\in[0,T]$. In order to guarantee that (H1) and (H2) are fulfilled we assume that the coefficients $\lambda_i$ are real valued for $i=1,2$ and we regularise the coefficients $\ell_{ii}$, $i=1,2$ via a logarithmic scale, i.e., $\omega^{-1}(\eps)=\ln(\eps^{-1})$. It is easy to solve the corresponding Eikonal equation and compute the phase functions $\phi_{i,\eps}$. Indeed,
	\[
	\phi_{i,\eps}(t,x,\xi)=x\xi+\xi\int_0^t \lambda_{i,\eps}(s)\, ds.
	\]
	It follows that 
	\begin{align*}
	|\partial_{\xi} \phi_{i, \eps}(t, x,\xi)-\partial_{\xi}\phi_{i, \eps}(t, y,\xi)| &= |x-y| \,\,\,\, \text{ for } \,\, x,y\in \mathbb{R}^n, \,\, \xi \in \mathbb{R}^n,  \\
	|\partial_y\phi_{i, \eps}(t, y,\xi)-\partial_y \phi_{i, \eps}(t, y,\eta)| &= |\xi - \eta|\,\,\,\, \text{ for } \,\,  y\in \mathbb{R}^n, \,\, \xi, \eta \in \mathbb{R}^n, \nonumber
	\end{align*}
	and therefore the net $C(\eps)$ required in (H3) is identically equal to $1$. In addition, since the regularisation is not effecting the variable $x$ all the constants involved in (H3) are independent of the parameter $\eps$ are therefore (H3) is trivially satisfied. It follows that our Cauchy problem is well-posed in the very weak sense.

	
	\section{Appendix: $L^2$-boundedness of Fourier integral operators}
	
	In this appendix we state an important  theorem by Ruzhansky-Sugimoto \cite{RuSu:06} on the global $L^2$ boundedness of FIO for amplitudes which are independent of the variable $x$. This theorem will be later applied to integrated FIOs.
	\begin{theorem}[Ruzhansky-Sugimoto \cite{RuSu:06}]\label{thm-L2-bound}
		Let operator $T$ be defined by 
		\begin{equation}
		Tu(x) = \int_{\mathbb{R}^n} \int_{\mathbb{R}^n} e^{i(x\cdot \xi + \phi(y,\xi))} A(y,\xi) u(y) dyd\xi,
		\end{equation}
		where $A(y,\xi) \in C^{\infty}(\mathbb{R}^n_{y}\times \mathbb{R}^n_{\xi})$, and $\phi(y,\xi)\in C^{\infty}(\mathbb{R}^n_y\times \mathbb{R}^n_{\xi})$ is a real-valued function. Assume that $|\partial_y^{\alpha}\partial_{\xi}^{\beta}A(y,\xi)|$ is uniformly bounded
		for $|\alpha|, |\beta|\leq 2n+1$. Also assume that
		\begin{align}\label{eq-22}
		|\partial_{\xi} \phi(x,\xi)-\partial_{\xi}\phi(y,\xi)| &\geq C|x-y| \,\,\,\, \text{ for } \,\, x,y\in \mathbb{R}^n, \,\, \xi \in \mathbb{R}^n, \\
		|\partial_y\phi(y,\xi)-\partial_y \phi(y,\eta)| &\geq C|\xi - \eta| \,\,\,\, \text{ for } \,\,  y\in \mathbb{R}^n, \,\, \xi, \eta \in \mathbb{R}^n, \nonumber
		\end{align}
		and that 
		\begin{equation}\label{eq-23}
		|\partial_y^{\alpha}\partial_{\xi}\phi(y,\xi)| \leq C_{\alpha}, \,\, 	|\partial_y\partial_{\xi}^{\beta}\phi(y,\xi)| \leq C_{\beta},
		\end{equation}
		for $1\leq |\alpha|, |\beta|\leq 2n+2$. Then the operator $T$ is $L^2(\mathbb{R}^n)$-bounded, and satisfies
		\begin{equation*}
		|| T||_{L^2 \rightarrow L^2} \leq C'\sup_{|\alpha|,|\beta|\leq 2n+1} || \partial_y^{\alpha} \partial_{\xi}^{\beta}A(y,\xi)||_{L^{\infty}(\mathbb{R}^n_y\times \mathbb{R}^n_{\xi})},
		\end{equation*}
		where  $C'$ depends on the constants in \eqref{eq-22} and \eqref{eq-23}.
	\end{theorem}
	
	Note that since $\Gamma_x=\Gamma_y=\mathbb{R}^n$ in Theorem \ref{thm-L2-bound}, we may write the following assumption 
	\begin{equation*}
	|\det \partial_x \partial_y \phi(x,y)|\geq C >0, \,\,\, \text{ for all } x,y \in \mathbb{R}^n,
	\end{equation*}
	instead of \eqref{eq-22}.
	
	Let us now consider operators of the type  
	\begin{align*}
	G u(t,x) &= \int_0^t \int_{\mathbb{R}^n} e^{i\phi(t,x,\xi)}A(t,s,x) \hat{u}(s,\xi)d\xi ds \\
	& = \int_0^t \int_{\mathbb{R}^n}\int_{\mathbb{R}^n} e^{i\phi(t,x,\xi)- iy\cdot \xi}A(t,s,x) u(s,y)dyd\xi ds\\
	& = \int_0^t Ku(t,s,x)ds,
	\end{align*}
	where 
	\begin{align*}
	Ku(s,x) = \int_{\mathbb{R}^n}\int_{\mathbb{R}^n} e^{i\phi(t,x,\xi)- iy\cdot \xi}A(t,s,x) u(s,y)dyd\xi. 
	\end{align*}
	It follows that
	\begin{align*}
	|| Gu(t)||_{L^2(\mathbb{R}^n_x)} = \int_0^t || K g(s)||_{L^2(\mathbb{R}^n_x)} ds &\leq t ||K||_{L^2 \rightarrow L^2} ||u||_{L^{\infty}(\mathbb{R}_s)L^2(\mathbb{R}^n_x)} \\
	& = t ||K^*||_{L^2 \rightarrow L^2} ||u||_{L^{\infty}(\mathbb{R}_s)L^2(\mathbb{R}^n_x)}. 
	\end{align*}
	In order to apply Theorem \ref{thm-L2-bound}, we compute the adjoint of operator $K$
	\begin{align*}
	(Ku,v) &= \int_{\mathbb{R}^n} \int_{\mathbb{R}^n}\int_{\mathbb{R}^n} e^{i\phi(t,x,\xi)- iy\cdot \xi}A(t,s,x) u(s,y)dyd\xi \overline{v}(s,x)dx \\ 
	&= \int_{\mathbb{R}^n} \overline{\int_{\mathbb{R}^n}\int_{\mathbb{R}^n} e^{i\phi(t,x,\xi)- iy\cdot \xi}A(t,s,x) \overline{v}(s,x)dxd\xi} u(s,y)dy \\
	&= \int_{\mathbb{R}^n} \int_{\mathbb{R}^n}\int_{\mathbb{R}^n} e^{-i\phi(t,x,\xi)+ iy\cdot \xi}A(t,s,x)v(s,x)dxd\xi u(s,y)dy \\
	& = (u,K^*v),
	\end{align*}
	where
	\begin{align*}
	K^*v(t,s,x) = \int_{\mathbb{R}^n}\int_{\mathbb{R}^n} e^{ ix\cdot \xi -i\phi(t,y,\xi)}A(t,s,y)v(s,y)dyd\xi.
	\end{align*}
	Now the adjoint operator $K^*$ has the same form as the operator $T$ in Theorem \ref{thm-L2-bound}.
	\begin{proposition}\label{Prop-L2}
		Let operator $K^*$ be defined by 
		\begin{equation}
		K^*u(t,s,x) = \int_{\mathbb{R}^n}\int_{\mathbb{R}^n} e^{ ix\cdot \xi -i\phi(t,y,\xi)}A(t,s,y)u(s,y)dyd\xi,
		\end{equation}
		where $A(t,s,y) \in C^{\infty}( \mathbb{R}^n_{y})$ is continuous with respect to $t,s\in[0,T]$ and $\phi(t, y,\xi)\in C^{\infty}(\mathbb{R}^n_y\times \mathbb{R}^n_{\xi})$ is a real-valued function continuous with respect to $t\in[0,T]$. Assume that $|\partial_y^{\gamma}A(s,t,y)|$ is uniformly bounded for $|\gamma|\leq 2n+1$. Also assume that
		\begin{align}
		|\partial_{\xi} \phi(t, x,\xi)-\partial_{\xi}\phi(t, y,\xi)| &\geq C|x-y| \,\,\,\, \text{ for } \,\, t\in[0,T],\, x,y\in \mathbb{R}^n, \,\, \xi \in \mathbb{R}^n, \label{eq-221} \\
		|\partial_y\phi(t, y,\xi)-\partial_y \phi(t, y,\eta)| &\geq C |\xi - \eta| \,\,\,\, \text{ for } \,\,  t\in[0,T],\, y\in \mathbb{R}^n, \,\, \xi, \eta \in \mathbb{R}^n, \nonumber
		\end{align}
		and that 
		\begin{equation}\label{eq-231}
		|\partial_y^{\alpha}\partial_{\xi}\phi(t, y,\xi)| \leq C_{\alpha}, \,\, 	|\partial_y\partial_{\xi}^{\beta}\phi(t, y,\xi)| \leq C_{\beta},
		\end{equation}
		for $1\leq |\alpha|, |\beta|\leq 2n+2$, uniformly with respect to all the variables. Then the operator $K^*(t,s,\cdot)$ is $L^2(\mathbb{R}^n)$-bounded, and satisfies
		\begin{equation*}
		\sup_{t,s\in[0,T]}|| K^*(t,s,\cdot)||_{L^2 \rightarrow L^2} \leq C'\sup_{t,s\in[0,T]}\sup_{|\gamma|\leq 2n+1} || \partial_y^{\gamma} A(t,s,y)||_{L^{\infty}(\mathbb{R}^n_y)},
		\end{equation*}
		where  $C'$ depends on the constants in \eqref{eq-221} and \eqref{eq-231}.
	\end{proposition}

\end{document}